\documentclass[1p]{elsarticle}
\usepackage{amsthm,amstext, amsmath,latexsym,amsbsy,amssymb,amscd}
\usepackage{graphicx}

\newtheorem{remark}{Remark}[section]
\newtheorem{lemma}{Lemma}[section] 
\newtheorem{example}{Example}[section]
\numberwithin{equation}{section}
\newtheorem{theorem}{Theorem}[section] 

\newcommand{\h}{\hspace*{.24in}}

\allowdisplaybreaks
\journal{}

\begin{document}

\begin{frontmatter}

\title{Convergence Properties of Overlapping Schwarz Domain Decomposition Algorithms}
\author{Minh-Binh TRAN}
\address{Laboratoire Analyse G\'eom\'etrie et Applications\\
Institut Galil\'ee, Universit\'e Paris 13, France
}
\ead{binh@math.univ-paris13.fr}
\begin{abstract} 
In this paper, we partially answer open questions about the convergence of overlapping Schwarz methods. We prove that overlapping Schwarz methods with Dirichlet transmission conditions for semilinear elliptic and parabolic equations always converge. While overlapping Schwarz methods with Robin transmission conditions only converge for semilinear parabolic equations, but not for semilinear elliptic ones. We then provide some conditions so that overlapping Schwarz methods with Robin transmission conditions converge for semilinear elliptic equations. Our new techniques can also be potentially applied to others kinds of partial differential equations.
\end{abstract}
\end{frontmatter}
\section{Introduction}
\h The Schwarz domain decomposition methods are procedures to solve partial differential equations in parallel, where each iteration involves the solutions of the original equation on smaller subdomains. The alternate method was originally proposed by H. A. Schwarz \cite{Schwarz:1890:GMA} in 1870 as a technique to prove the existence of a solution to the Laplace equation on a domain which is a combination of a rectangle and a circle. The idea was then used and extended by P. L. Lions \cite{Lions:1987:OSA}, \cite{Lions:1989:OSA}, \cite{Lions:1990:OSA} to parallel algorithms for solving partial differential equations. Since then, many kind of domain decomposition methods have been developed, to improve the performance of the classical domain decomposition method. However, the convergence for domain decomposition methods still remains an open question.
\\\h Many techniques have been developed to prove the convergence of classical Schwarz methods, or Schwarz methods with Dirichlet transmission conditions. One of the first techniques, used by P. L. Lions in \cite{Lions:1987:OSA}, is the iterated projections for linear Laplace equation and linear Stoke equation. The idea is to prove that classical Schwarz methods for these equations are equivalent to sequences of projections in Hilbert Spaces. In the same paper, P. L. Lions also showed that the Schwarz sequences for nonlinear monotone elliptic equations are related to classical minimization methods over product spaces and proposed to use Schwarz methods for evolution equations. This idea was then used by L. Badea in \cite{Badea:1991:OSA} to prove the convergence of classical Schwarz methods for nonlinear monotone elliptic problems.
\\ Following the pinoneering work of P. L. Lions, in the papers \cite{Gander:1999:WRA}, \cite{Giladi:2002:STD}, \cite{GanderStuart:1998:STC}, E. Giladi, H. B. Keller, A. Stuart and M. Gander used Fourier and Laplace transforms, together with some explicit calculation to study classical Schwaz methods for some 1-dimensional evolution equations, with constant coefficients. Later, by using a maximum principle argument, M. Gander and H. Zhao proved that classical Schwarz method converges for the n-dimensional linear heat equation \cite{Gander:2002:OSW}.
\\ Another technique to study the convergence of classical Schwarz methods is to use the idea of upper-lower solutions methods, with initial guess to be upper or lower solutions of the equations. This special class of domain decomposition methods with monotone iterations has been studied by S. H. Lui in \cite{Lui:2002:OLM}, \cite{Lui:2003:OMI}, \cite{Lui:2001:OMS}. Although many techniques have been developed to study the convergence problem of classical Schwarz methods, the problem with nonlinear equations in n-dimension and general multi-subdomains is still open.
\\\h A new class of Schwarz algorithms, in which Dirichlet transmission condition is replaced by Robin ones, has been studied recently in order to improve the performance of classical methods. The new algorithms are called optimized Schwarz methods since there are some parameters we can optimize to get faster algorithms. In 1989, P.L.
Lions (see \cite{Lions:1989:OSA}, \cite{Lions:1990:OSA}) established the convergence of nonoverlapping optimized Schwarz methods with Robin transmission conditions by using an energy argument. Later, J. D. Benamou and B. Depres in \cite{Benamou:1997:DDM} used this technique to study the convergence of nonoverlapping  optimized Schwarz methods for Helmholtz equation. Energy estimates have then become a very powerful technique to prove the convergence of nonoverlapping optimized Schwarz methods with Robin transmission conditions (see \cite{HalpernSzeftel:2008:OQO}).
\\ However, the convergence problem of overlapping optimized Schwarz methods, even for linear problems, still remains an open problem up to now. J.-H. Kimn \cite{Kimn:2005:CTO}, proved the convergence of an overlapping optimized Schwarz method for Poisson equation with Robin boundary data, 
\begin{equation*}
\left \{ \begin{array}{ll}-\Delta u=f \mbox{ in }\Omega,\vspace{.1in}\\ 
\frac{\partial u}{\partial n}+pu=g \mbox{ on } \partial\Omega.\end{array}\right. 
\end{equation*}
He proved that there is an $p_0 > 0$ such that the Schwarz iterations with Robin transmissions conditions converge for any Robin parameter $0 <p <p_0$. In \cite{LoiselSzyld:2010:OGC}, S. Loisel and D. B. Szyld extended the technique of J.-H. Kimn for the following equation
\begin{equation*}
\left \{ \begin{array}{ll}-\nabla(a\nabla u) +c u=f \mbox{ in }\Omega,\vspace{.1in}\\ 
u=0 \mbox{ on } \partial\Omega,\end{array}\right. 
\end{equation*}
where $a$ is a $C^1$-function and $c$ is positive and belongs to $L^{\infty}(\Omega)$. The same constant $p$ is kept for all transmission operators and some conditions on the boundaries of the subdomains are then imposed. 
\\ A proof of convergence based on semi-classical analysis for overlapping optimized Schwarz methods with rectangle subdomains, linear advection diffusion equations on the half plane was given in \cite{Nataf:1998:CDD}. 
\\ Another technique is to use Fourier transform. This technique cannot be used to study the convergence of Schwarz methods for nonlinear problems and for general subdomains, but convergence rates can be obtained. Changing the boundary conditions will change the values of the convergence rates and then improve the performance of the algorithms, which proposes a new problem: the problem of optimizing the convergence rates. In \cite{JaphetNataf:2001:TBI}, \cite{Bennequin:2009:AHB}, \cite{Gander:2007:OSW} the authors showed that the problem of optimizing the convergence rates is in fact a new class of best approximation problems and suggested a new method to solve it.
\\\h In this paper, we present convergence proofs  of overlapping classical and optimized Schwarz methods for elliptic and parabolic semilinear equations, in general forms, for general multi-subdomains. We prove that Schwarz methods with Dirichlet and Robin transmission conditions always converge for parabolic equations; since with parabolic equations the time variable can be controlled easily. However, Schwarz methods with Robin transmission conditions do not converge for elliptic equations, while classical Schwarz methods always converge. We can see from Remark 3.1 that given a Schwarz algorithm with a specific Robin transmission condition, there exists a class of elliptic equations where the algorithm is unstable. A condition of convergence is then supplied: Schwarz methods with Robin transmission conditions for elliptic equations will converge if we multiply Robin parameters by a number large enough, and this can also be seen from Example 3.1. The techniques used in our proofs can also be used to prove the convergence of Schwarz methods for many other kinds of partial differential equations.
\\\h The paper is organized as follows.
\begin{itemize}
\item Section $2$ is devoted to the convergence properties of Schwarz methods for semilinear parabolic equations. Section $2.1$ gives the definition of the Schwarz algorithms for semilinear parabolic equations, and states the two theorems of convergence. Theorem $2.1$ announces that classical Schwarz algorithms always converge with semilinear parabolic equations and its proof can be found in section $2.2$. Theorem $2.2$ is about the convergence of Schwarz algorithms with Robin transmission conditions and the proof is then given later in section $2.3$.
\item In section $3$, we discuss the convergence properties of Schwarz methods for semilinear elliptic equations. Definitions of the algorithms, the two convergence theorems $3.1$, $3.2$, together with the counterexample $3.1$ is announced in section $3.1$. Section $3.2$ and $3.3$ contain the proofs of the two theorems.
\end{itemize}
\section{Convergence for Semilinear Parabolic Equations}
\h We introduce the abbreviation $\partial_{i,j}=\frac{\partial^2}{\partial x_i\partial x_j}$, $\partial_t=\frac{\partial}{\partial t}$ and $\partial_i=\frac{\partial }{\partial x_i}$ and consider a general semilinear parabolic equation
\begin{equation}
\label{2e1}
\left \{ \begin{array}{ll}{\partial_t u}-\sum_{i,j=1}^n\partial_j(a_{i,j}{\partial_{i} u})+\sum_{i=1}^nb_{i}{\partial_i u}+cu=F(x,t,u)\mbox{ in } \Omega\times(0,\infty),\vspace{.1in}\\ 
u(x,t)=g(x,t) \mbox{ on } \partial\Omega\times(0,\infty),\vspace{.1in}\\ 
u(x,0)=g(x,0) \mbox{ on } \Omega,\end{array}\right. 
\end{equation}
where $\Omega$ is a bounded and smooth domain in  $\mathbb{R}^{n}$. The coefficients $a_{i,j}$, $b_i$, $c$ are functions of the space variable $x$, with the following properties
\\ (A1) The functions $a_{i,j}$, $b_i$, $c$ are in $C^2(\mathbb{R}^{n})$. 
\\ (A2) For all $i,j$ in $\{1,\dots,I\}$, $a_{i,j}(x)=a_{j,i}(x)$. There exist strictly positive numbers $\lambda$, $\Lambda$ such that $A=(a_{i,j}(x))\geq\lambda I$ in the sense of symmetric positive definite matrices and $|a_{i,j}(x)|<\Lambda$  in $\Omega$.
\\ (A3) $g$ is in $C^{2}(\mathbb{R}^{n+1})$ and and $F$ is uniformly Lipschitz in the third variable, i.e. there exists $C>0$, such that \\\h $\forall$ $t$ $\in$ $\mathbb{R}$, $\forall$ $x$ $\in$ $\mathbb{R}^n$, $|F(x,t,z)-F(x,t,z')|\leq C|z-z'|$, $\forall$ $z$, $z'$ $\in$ $\mathbb{R}$.
\\ With Conditions $(A1)$, $(A2)$ and $(A3)$, Equation $(\ref{2e1})$ has a unique bounded solution $u$ in $C^{2,1}(\overline{\Omega\times(0,\infty)})$, i.e $\partial_{i,j}\partial_tu$ belongs to $C(\overline{\Omega\times(0,\infty)})$ for all $i,j$ in $\{1,\dots,n\}$. The proof of this result can be found in some classical books like \cite{Friedman:1964:PDE}, \cite{Lieberman:1996:SOP}.
\\\h The domain $\Omega$ is divided into $I$ smooth overlapping subdomains $\{\Omega_l\}_{l\in\{1,I\}}$, such that
$$\cup_{l=1}^{n}\Omega_l=\Omega;$$ $$(\partial\Omega_l\backslash\partial\Omega)\cap(\partial\Omega_{l'}\backslash\partial\Omega)=\O,~~\forall~~l,l'\in \{1,\dots,I\}, ~~l\ne l';$$ 
and
$$\forall l\in \{1,\dots,I\}, \forall l',l''\in J_l,l''\ne l',~~~\Omega_{l'}\cap\Omega_{l''}=\O,$$
where
$$J_l=\{l'|\Omega_{l'}\cap\Omega_l\ne\O\}.$$
For any $l$ in $J$, for $l'\in J_l$, $\Gamma_{l,l'}$ is the set $(\partial\Omega_l\backslash\partial\Omega)\cap\overline\Omega_{l'}$. 
\begin{remark}
\begin{figure}[!ht]
\centering
%WinTpicVersion3.08
\unitlength 1pt
\begin{picture}(339.2991,205.5712)(128.0374,-275.2804)
% BOX 3 0 3 0
% 2 1800 980 6570 3870
% 
\special{pn 4}%
\special{pa 1772 965}%
\special{pa 6467 965}%
\special{pa 6467 3810}%
\special{pa 1772 3810}%
\special{pa 1772 965}%
\special{fp}%
% BOX 3 0 3 0
% 2 1810 980 1810 990
% 
\special{pn 4}%
\special{pa 1782 965}%
\special{pa 1782 965}%
\special{pa 1782 975}%
\special{pa 1782 975}%
\special{pa 1782 965}%
\special{fp}%
% BOX 0 1 3 0
% 2 1820 990 4920 3850
% 
\special{pn 20}%
\special{pa 1792 975}%
\special{pa 4843 975}%
\special{pa 4843 3790}%
\special{pa 1792 3790}%
\special{pa 1792 975}%
\special{da 0.070}%
% CIRCLE 1 2 3 0
% 4 2820 2480 3180 3130 3180 3130 3180 3120
% 
\special{pn 13}%
\special{ar 2776 2441 732 732  1.0584069 1.0650050}%
% CIRCLE 3 0 1 0
% 4 2660 2120 2990 2750 2990 2750 2990 2750
% 
\special{pn 4}%
\special{sh 0.300}%
\special{ar 2619 2087 700 700  0.0000000 6.2831853}%
% CIRCLE 0 1 3 0
% 4 2610 2020 2690 2370 2690 2370 2690 2370
% 
\special{pn 20}%
\special{ar 2569 1989 354 354  0.0000000 0.1671309}%
\special{ar 2569 1989 354 354  0.2674095 0.4345404}%
\special{ar 2569 1989 354 354  0.5348189 0.7019499}%
\special{ar 2569 1989 354 354  0.8022284 0.9693593}%
\special{ar 2569 1989 354 354  1.0696379 1.2367688}%
\special{ar 2569 1989 354 354  1.3370474 1.5041783}%
\special{ar 2569 1989 354 354  1.6044568 1.7715877}%
\special{ar 2569 1989 354 354  1.8718663 2.0389972}%
\special{ar 2569 1989 354 354  2.1392758 2.3064067}%
\special{ar 2569 1989 354 354  2.4066852 2.5738162}%
\special{ar 2569 1989 354 354  2.6740947 2.8412256}%
\special{ar 2569 1989 354 354  2.9415042 3.1086351}%
\special{ar 2569 1989 354 354  3.2089136 3.3760446}%
\special{ar 2569 1989 354 354  3.4763231 3.6434540}%
\special{ar 2569 1989 354 354  3.7437326 3.9108635}%
\special{ar 2569 1989 354 354  4.0111421 4.1782730}%
\special{ar 2569 1989 354 354  4.2785515 4.4456825}%
\special{ar 2569 1989 354 354  4.5459610 4.7130919}%
\special{ar 2569 1989 354 354  4.8133705 4.9805014}%
\special{ar 2569 1989 354 354  5.0807799 5.2479109}%
\special{ar 2569 1989 354 354  5.3481894 5.5153203}%
\special{ar 2569 1989 354 354  5.6155989 5.7827298}%
\special{ar 2569 1989 354 354  5.8830084 6.0501393}%
\special{ar 2569 1989 354 354  6.1504178 6.2832853}%
% BOX 3 1 3 0
% 2 4570 990 6560 3850
% 
\special{pn 4}%
\special{pa 4499 975}%
\special{pa 6457 975}%
\special{pa 6457 3790}%
\special{pa 4499 3790}%
\special{pa 4499 975}%
\special{da 0.070}%
\end{picture}%

\caption{A good way of dividing $\Omega$}
\end{figure}
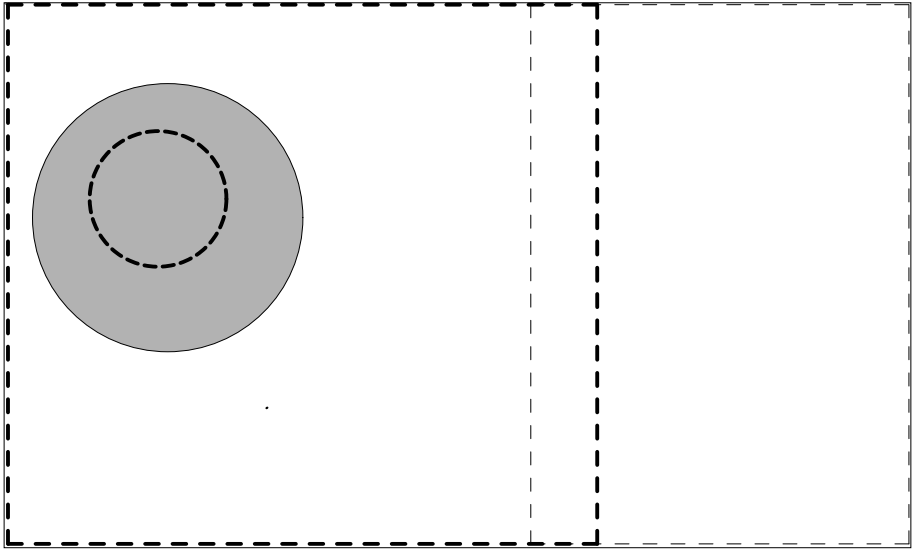
\begin{figure}[!ht]
\centering
\unitlength 0.1in
\begin{picture}( 38.5000, 26.0000)( 25.7000,-38.0000)
% BOX 2 5 3 0
% 2 2580 1200 6410 3790
% 
\special{pn 8}%
\special{pa 2580 1200}%
\special{pa 6410 1200}%
\special{pa 6410 3790}%
\special{pa 2580 3790}%
\special{pa 2580 1200}%
\special{ip}%
% BOX 1 0 3 0
% 2 2570 2830 6420 3800
% 
\special{pn 13}%
\special{pa 2570 2830}%
\special{pa 6420 2830}%
\special{pa 6420 3800}%
\special{pa 2570 3800}%
\special{pa 2570 2830}%
\special{fp}%
% LINE 2 0 3 0
% 60 4010 2830 3650 3190 3950 2830 3590 3190 3890 2830 3530 3190 3830 2830 3470 3190 3770 2830 3410 3190 3710 2830 3350 3190 3650 2830 3290 3190 3590 2830 3230 3190 3530 2830 3170 3190 3470 2830 3110 3190 3410 2830 3050 3190 3350 2830 2990 3190 3290 2830 2930 3190 3230 2830 2870 3190 3170 2830 2810 3190 3110 2830 2750 3190 3050 2830 2690 3190 2990 2830 2630 3190 2930 2830 2580 3180 2630 3190 2580 3240 2690 3190 2580 3300 2750 3190 2580 3360 2810 3190 2580 3420 2870 3190 2580 3480 2930 3190 2580 3540 2990 3190 2580 3600 3050 3190 2580 3660 3110 3190 2580 3720 3170 3190 2580 3780 3230 3190 2630 3790
% 
\special{pn 8}%
\special{pa 4010 2830}%
\special{pa 3650 3190}%
\special{fp}%
\special{pa 3950 2830}%
\special{pa 3590 3190}%
\special{fp}%
\special{pa 3890 2830}%
\special{pa 3530 3190}%
\special{fp}%
\special{pa 3830 2830}%
\special{pa 3470 3190}%
\special{fp}%
\special{pa 3770 2830}%
\special{pa 3410 3190}%
\special{fp}%
\special{pa 3710 2830}%
\special{pa 3350 3190}%
\special{fp}%
\special{pa 3650 2830}%
\special{pa 3290 3190}%
\special{fp}%
\special{pa 3590 2830}%
\special{pa 3230 3190}%
\special{fp}%
\special{pa 3530 2830}%
\special{pa 3170 3190}%
\special{fp}%
\special{pa 3470 2830}%
\special{pa 3110 3190}%
\special{fp}%
\special{pa 3410 2830}%
\special{pa 3050 3190}%
\special{fp}%
\special{pa 3350 2830}%
\special{pa 2990 3190}%
\special{fp}%
\special{pa 3290 2830}%
\special{pa 2930 3190}%
\special{fp}%
\special{pa 3230 2830}%
\special{pa 2870 3190}%
\special{fp}%
\special{pa 3170 2830}%
\special{pa 2810 3190}%
\special{fp}%
\special{pa 3110 2830}%
\special{pa 2750 3190}%
\special{fp}%
\special{pa 3050 2830}%
\special{pa 2690 3190}%
\special{fp}%
\special{pa 2990 2830}%
\special{pa 2630 3190}%
\special{fp}%
\special{pa 2930 2830}%
\special{pa 2580 3180}%
\special{fp}%
\special{pa 2630 3190}%
\special{pa 2580 3240}%
\special{fp}%
\special{pa 2690 3190}%
\special{pa 2580 3300}%
\special{fp}%
\special{pa 2750 3190}%
\special{pa 2580 3360}%
\special{fp}%
\special{pa 2810 3190}%
\special{pa 2580 3420}%
\special{fp}%
\special{pa 2870 3190}%
\special{pa 2580 3480}%
\special{fp}%
\special{pa 2930 3190}%
\special{pa 2580 3540}%
\special{fp}%
\special{pa 2990 3190}%
\special{pa 2580 3600}%
\special{fp}%
\special{pa 3050 3190}%
\special{pa 2580 3660}%
\special{fp}%
\special{pa 3110 3190}%
\special{pa 2580 3720}%
\special{fp}%
\special{pa 3170 3190}%
\special{pa 2580 3780}%
\special{fp}%
\special{pa 3230 3190}%
\special{pa 2630 3790}%
\special{fp}%
% LINE 2 0 3 1
% 60 3290 3190 2690 3790 3350 3190 2750 3790 3410 3190 2810 3790 3470 3190 2870 3790 3530 3190 2930 3790 3590 3190 2990 3790 3650 3190 3050 3790 3710 3190 3110 3790 3770 3190 3170 3790 3830 3190 3230 3790 3890 3190 3290 3790 3950 3190 3350 3790 4010 3190 3410 3790 4070 3190 3470 3790 4130 3190 3530 3790 4190 3190 3590 3790 4250 3190 3650 3790 4310 3190 3710 3790 4370 3190 3770 3790 4430 3190 3830 3790 4490 3190 3890 3790 4550 3190 3950 3790 4590 3210 4010 3790 4590 3270 4070 3790 4590 3330 4130 3790 4590 3390 4190 3790 4590 3450 4250 3790 4590 3510 4310 3790 4610 3550 4370 3790 4670 3550 4430 3790
% 
\special{pn 8}%
\special{pa 3290 3190}%
\special{pa 2690 3790}%
\special{fp}%
\special{pa 3350 3190}%
\special{pa 2750 3790}%
\special{fp}%
\special{pa 3410 3190}%
\special{pa 2810 3790}%
\special{fp}%
\special{pa 3470 3190}%
\special{pa 2870 3790}%
\special{fp}%
\special{pa 3530 3190}%
\special{pa 2930 3790}%
\special{fp}%
\special{pa 3590 3190}%
\special{pa 2990 3790}%
\special{fp}%
\special{pa 3650 3190}%
\special{pa 3050 3790}%
\special{fp}%
\special{pa 3710 3190}%
\special{pa 3110 3790}%
\special{fp}%
\special{pa 3770 3190}%
\special{pa 3170 3790}%
\special{fp}%
\special{pa 3830 3190}%
\special{pa 3230 3790}%
\special{fp}%
\special{pa 3890 3190}%
\special{pa 3290 3790}%
\special{fp}%
\special{pa 3950 3190}%
\special{pa 3350 3790}%
\special{fp}%
\special{pa 4010 3190}%
\special{pa 3410 3790}%
\special{fp}%
\special{pa 4070 3190}%
\special{pa 3470 3790}%
\special{fp}%
\special{pa 4130 3190}%
\special{pa 3530 3790}%
\special{fp}%
\special{pa 4190 3190}%
\special{pa 3590 3790}%
\special{fp}%
\special{pa 4250 3190}%
\special{pa 3650 3790}%
\special{fp}%
\special{pa 4310 3190}%
\special{pa 3710 3790}%
\special{fp}%
\special{pa 4370 3190}%
\special{pa 3770 3790}%
\special{fp}%
\special{pa 4430 3190}%
\special{pa 3830 3790}%
\special{fp}%
\special{pa 4490 3190}%
\special{pa 3890 3790}%
\special{fp}%
\special{pa 4550 3190}%
\special{pa 3950 3790}%
\special{fp}%
\special{pa 4590 3210}%
\special{pa 4010 3790}%
\special{fp}%
\special{pa 4590 3270}%
\special{pa 4070 3790}%
\special{fp}%
\special{pa 4590 3330}%
\special{pa 4130 3790}%
\special{fp}%
\special{pa 4590 3390}%
\special{pa 4190 3790}%
\special{fp}%
\special{pa 4590 3450}%
\special{pa 4250 3790}%
\special{fp}%
\special{pa 4590 3510}%
\special{pa 4310 3790}%
\special{fp}%
\special{pa 4610 3550}%
\special{pa 4370 3790}%
\special{fp}%
\special{pa 4670 3550}%
\special{pa 4430 3790}%
\special{fp}%
% LINE 2 0 3 2
% 60 4730 3550 4490 3790 4790 3550 4550 3790 4850 3550 4610 3790 4910 3550 4670 3790 4970 3550 4730 3790 5030 3550 4790 3790 5080 3560 4850 3790 5140 3560 4910 3790 5200 3560 4970 3790 5260 3560 5030 3790 5320 3560 5090 3790 5380 3560 5150 3790 5440 3560 5210 3790 5500 3560 5270 3790 5560 3560 5330 3790 5620 3560 5390 3790 5680 3560 5450 3790 5740 3560 5510 3790 5800 3560 5570 3790 5860 3560 5630 3790 5920 3560 5690 3790 5980 3560 5750 3790 6030 3570 5810 3790 6090 3570 5870 3790 6150 3570 5930 3790 6210 3570 5990 3790 6270 3570 6050 3790 6330 3570 6110 3790 6390 3570 6170 3790 6410 3610 6230 3790
% 
\special{pn 8}%
\special{pa 4730 3550}%
\special{pa 4490 3790}%
\special{fp}%
\special{pa 4790 3550}%
\special{pa 4550 3790}%
\special{fp}%
\special{pa 4850 3550}%
\special{pa 4610 3790}%
\special{fp}%
\special{pa 4910 3550}%
\special{pa 4670 3790}%
\special{fp}%
\special{pa 4970 3550}%
\special{pa 4730 3790}%
\special{fp}%
\special{pa 5030 3550}%
\special{pa 4790 3790}%
\special{fp}%
\special{pa 5080 3560}%
\special{pa 4850 3790}%
\special{fp}%
\special{pa 5140 3560}%
\special{pa 4910 3790}%
\special{fp}%
\special{pa 5200 3560}%
\special{pa 4970 3790}%
\special{fp}%
\special{pa 5260 3560}%
\special{pa 5030 3790}%
\special{fp}%
\special{pa 5320 3560}%
\special{pa 5090 3790}%
\special{fp}%
\special{pa 5380 3560}%
\special{pa 5150 3790}%
\special{fp}%
\special{pa 5440 3560}%
\special{pa 5210 3790}%
\special{fp}%
\special{pa 5500 3560}%
\special{pa 5270 3790}%
\special{fp}%
\special{pa 5560 3560}%
\special{pa 5330 3790}%
\special{fp}%
\special{pa 5620 3560}%
\special{pa 5390 3790}%
\special{fp}%
\special{pa 5680 3560}%
\special{pa 5450 3790}%
\special{fp}%
\special{pa 5740 3560}%
\special{pa 5510 3790}%
\special{fp}%
\special{pa 5800 3560}%
\special{pa 5570 3790}%
\special{fp}%
\special{pa 5860 3560}%
\special{pa 5630 3790}%
\special{fp}%
\special{pa 5920 3560}%
\special{pa 5690 3790}%
\special{fp}%
\special{pa 5980 3560}%
\special{pa 5750 3790}%
\special{fp}%
\special{pa 6030 3570}%
\special{pa 5810 3790}%
\special{fp}%
\special{pa 6090 3570}%
\special{pa 5870 3790}%
\special{fp}%
\special{pa 6150 3570}%
\special{pa 5930 3790}%
\special{fp}%
\special{pa 6210 3570}%
\special{pa 5990 3790}%
\special{fp}%
\special{pa 6270 3570}%
\special{pa 6050 3790}%
\special{fp}%
\special{pa 6330 3570}%
\special{pa 6110 3790}%
\special{fp}%
\special{pa 6390 3570}%
\special{pa 6170 3790}%
\special{fp}%
\special{pa 6410 3610}%
\special{pa 6230 3790}%
\special{fp}%
% LINE 2 0 3 3
% 60 6410 3670 6290 3790 6410 3730 6350 3790 4910 3190 4590 3510 4970 3190 4610 3550 5390 2830 4670 3550 5450 2830 4730 3550 5510 2830 4790 3550 5570 2830 4850 3550 5630 2830 4910 3550 5690 2830 4970 3550 5750 2830 5030 3550 5810 2830 5080 3560 5870 2830 5140 3560 5930 2830 5200 3560 5990 2830 5260 3560 6050 2830 5320 3560 6110 2830 5380 3560 6170 2830 5440 3560 6230 2830 5500 3560 6290 2830 5560 3560 6350 2830 5620 3560 6400 2840 5680 3560 6410 2890 5740 3560 6410 2950 5800 3560 6410 3010 5860 3560 6410 3070 5920 3560 6410 3130 5980 3560 6410 3190 6030 3570 6410 3250 6090 3570 6410 3310 6150 3570
% 
\special{pn 8}%
\special{pa 6410 3670}%
\special{pa 6290 3790}%
\special{fp}%
\special{pa 6410 3730}%
\special{pa 6350 3790}%
\special{fp}%
\special{pa 4910 3190}%
\special{pa 4590 3510}%
\special{fp}%
\special{pa 4970 3190}%
\special{pa 4610 3550}%
\special{fp}%
\special{pa 5390 2830}%
\special{pa 4670 3550}%
\special{fp}%
\special{pa 5450 2830}%
\special{pa 4730 3550}%
\special{fp}%
\special{pa 5510 2830}%
\special{pa 4790 3550}%
\special{fp}%
\special{pa 5570 2830}%
\special{pa 4850 3550}%
\special{fp}%
\special{pa 5630 2830}%
\special{pa 4910 3550}%
\special{fp}%
\special{pa 5690 2830}%
\special{pa 4970 3550}%
\special{fp}%
\special{pa 5750 2830}%
\special{pa 5030 3550}%
\special{fp}%
\special{pa 5810 2830}%
\special{pa 5080 3560}%
\special{fp}%
\special{pa 5870 2830}%
\special{pa 5140 3560}%
\special{fp}%
\special{pa 5930 2830}%
\special{pa 5200 3560}%
\special{fp}%
\special{pa 5990 2830}%
\special{pa 5260 3560}%
\special{fp}%
\special{pa 6050 2830}%
\special{pa 5320 3560}%
\special{fp}%
\special{pa 6110 2830}%
\special{pa 5380 3560}%
\special{fp}%
\special{pa 6170 2830}%
\special{pa 5440 3560}%
\special{fp}%
\special{pa 6230 2830}%
\special{pa 5500 3560}%
\special{fp}%
\special{pa 6290 2830}%
\special{pa 5560 3560}%
\special{fp}%
\special{pa 6350 2830}%
\special{pa 5620 3560}%
\special{fp}%
\special{pa 6400 2840}%
\special{pa 5680 3560}%
\special{fp}%
\special{pa 6410 2890}%
\special{pa 5740 3560}%
\special{fp}%
\special{pa 6410 2950}%
\special{pa 5800 3560}%
\special{fp}%
\special{pa 6410 3010}%
\special{pa 5860 3560}%
\special{fp}%
\special{pa 6410 3070}%
\special{pa 5920 3560}%
\special{fp}%
\special{pa 6410 3130}%
\special{pa 5980 3560}%
\special{fp}%
\special{pa 6410 3190}%
\special{pa 6030 3570}%
\special{fp}%
\special{pa 6410 3250}%
\special{pa 6090 3570}%
\special{fp}%
\special{pa 6410 3310}%
\special{pa 6150 3570}%
\special{fp}%
% LINE 2 0 3 4
% 60 6410 3370 6210 3570 6410 3430 6270 3570 6410 3490 6330 3570 5330 2830 5010 3150 5270 2830 5010 3090 5210 2830 5010 3030 5150 2830 5000 2980 5090 2830 5000 2920 5030 2830 5000 2860 4850 3190 4590 3450 4790 3190 4590 3390 4730 3190 4590 3330 4670 3190 4590 3270 2870 2830 2580 3120 2810 2830 2580 3060 2750 2830 2580 3000 2690 2830 2580 2940 2630 2830 2580 2880 4070 2830 3710 3190 4130 2830 3770 3190 4190 2830 3830 3190 4250 2830 3890 3190 4310 2830 3950 3190 4370 2830 4010 3190 4430 2830 4070 3190 4490 2830 4130 3190 4550 2830 4190 3190 4590 2850 4250 3190 4590 2910 4310 3190 4590 2970 4370 3190
% 
\special{pn 8}%
\special{pa 6410 3370}%
\special{pa 6210 3570}%
\special{fp}%
\special{pa 6410 3430}%
\special{pa 6270 3570}%
\special{fp}%
\special{pa 6410 3490}%
\special{pa 6330 3570}%
\special{fp}%
\special{pa 5330 2830}%
\special{pa 5010 3150}%
\special{fp}%
\special{pa 5270 2830}%
\special{pa 5010 3090}%
\special{fp}%
\special{pa 5210 2830}%
\special{pa 5010 3030}%
\special{fp}%
\special{pa 5150 2830}%
\special{pa 5000 2980}%
\special{fp}%
\special{pa 5090 2830}%
\special{pa 5000 2920}%
\special{fp}%
\special{pa 5030 2830}%
\special{pa 5000 2860}%
\special{fp}%
\special{pa 4850 3190}%
\special{pa 4590 3450}%
\special{fp}%
\special{pa 4790 3190}%
\special{pa 4590 3390}%
\special{fp}%
\special{pa 4730 3190}%
\special{pa 4590 3330}%
\special{fp}%
\special{pa 4670 3190}%
\special{pa 4590 3270}%
\special{fp}%
\special{pa 2870 2830}%
\special{pa 2580 3120}%
\special{fp}%
\special{pa 2810 2830}%
\special{pa 2580 3060}%
\special{fp}%
\special{pa 2750 2830}%
\special{pa 2580 3000}%
\special{fp}%
\special{pa 2690 2830}%
\special{pa 2580 2940}%
\special{fp}%
\special{pa 2630 2830}%
\special{pa 2580 2880}%
\special{fp}%
\special{pa 4070 2830}%
\special{pa 3710 3190}%
\special{fp}%
\special{pa 4130 2830}%
\special{pa 3770 3190}%
\special{fp}%
\special{pa 4190 2830}%
\special{pa 3830 3190}%
\special{fp}%
\special{pa 4250 2830}%
\special{pa 3890 3190}%
\special{fp}%
\special{pa 4310 2830}%
\special{pa 3950 3190}%
\special{fp}%
\special{pa 4370 2830}%
\special{pa 4010 3190}%
\special{fp}%
\special{pa 4430 2830}%
\special{pa 4070 3190}%
\special{fp}%
\special{pa 4490 2830}%
\special{pa 4130 3190}%
\special{fp}%
\special{pa 4550 2830}%
\special{pa 4190 3190}%
\special{fp}%
\special{pa 4590 2850}%
\special{pa 4250 3190}%
\special{fp}%
\special{pa 4590 2910}%
\special{pa 4310 3190}%
\special{fp}%
\special{pa 4590 2970}%
\special{pa 4370 3190}%
\special{fp}%
% LINE 2 0 3 5
% 6 4590 3030 4430 3190 4590 3090 4490 3190 4590 3150 4550 3190
% 
\special{pn 8}%
\special{pa 4590 3030}%
\special{pa 4430 3190}%
\special{fp}%
\special{pa 4590 3090}%
\special{pa 4490 3190}%
\special{fp}%
\special{pa 4590 3150}%
\special{pa 4550 3190}%
\special{fp}%
% LINE 3 0 3 0
% 24 4970 2830 4610 3190 4910 2830 4590 3150 4850 2830 4590 3090 4790 2830 4590 3030 4730 2830 4590 2970 4670 2830 4590 2910 5000 2860 4670 3190 5000 2920 4730 3190 5000 2980 4790 3190 5010 3030 4850 3190 5010 3090 4910 3190 5010 3150 4970 3190
% 
\special{pn 4}%
\special{pa 4970 2830}%
\special{pa 4610 3190}%
\special{fp}%
\special{pa 4910 2830}%
\special{pa 4590 3150}%
\special{fp}%
\special{pa 4850 2830}%
\special{pa 4590 3090}%
\special{fp}%
\special{pa 4790 2830}%
\special{pa 4590 3030}%
\special{fp}%
\special{pa 4730 2830}%
\special{pa 4590 2970}%
\special{fp}%
\special{pa 4670 2830}%
\special{pa 4590 2910}%
\special{fp}%
\special{pa 5000 2860}%
\special{pa 4670 3190}%
\special{fp}%
\special{pa 5000 2920}%
\special{pa 4730 3190}%
\special{fp}%
\special{pa 5000 2980}%
\special{pa 4790 3190}%
\special{fp}%
\special{pa 5010 3030}%
\special{pa 4850 3190}%
\special{fp}%
\special{pa 5010 3090}%
\special{pa 4910 3190}%
\special{fp}%
\special{pa 5010 3150}%
\special{pa 4970 3190}%
\special{fp}%
% BOX 0 1 3 0
% 2 4370 1210 6420 3370
% 
\special{pn 20}%
\special{pa 4370 1210}%
\special{pa 6420 1210}%
\special{pa 6420 3370}%
\special{pa 4370 3370}%
\special{pa 4370 1210}%
\special{da 0.070}%
% LINE 3 2 3 0
% 60 5000 1760 4450 1210 5000 1820 4390 1210 5000 1880 4370 1250 5000 1940 4370 1310 5000 2000 4370 1370 5000 2060 4370 1430 5000 2120 4370 1490 5000 2180 4370 1550 5000 2240 4370 1610 5000 2300 4370 1670 5000 2360 4370 1730 5000 2420 4370 1790 5000 2480 4370 1850 5000 2540 4370 1910 5000 2600 4370 1970 5000 2660 4370 2030 5000 2720 4370 2090 5000 2780 4370 2150 4990 2830 4370 2210 4930 2830 4370 2270 4870 2830 4370 2330 4810 2830 4370 2390 4750 2830 4370 2450 4690 2830 4370 2510 4630 2830 4370 2570 4570 2830 4370 2630 4510 2830 4370 2690 4450 2830 4370 2750 5050 2830 5000 2780 5110 2830 5000 2720
% 
\special{pn 4}%
\special{pa 5000 1760}%
\special{pa 4450 1210}%
\special{dt 0.027}%
\special{pa 5000 1820}%
\special{pa 4390 1210}%
\special{dt 0.027}%
\special{pa 5000 1880}%
\special{pa 4370 1250}%
\special{dt 0.027}%
\special{pa 5000 1940}%
\special{pa 4370 1310}%
\special{dt 0.027}%
\special{pa 5000 2000}%
\special{pa 4370 1370}%
\special{dt 0.027}%
\special{pa 5000 2060}%
\special{pa 4370 1430}%
\special{dt 0.027}%
\special{pa 5000 2120}%
\special{pa 4370 1490}%
\special{dt 0.027}%
\special{pa 5000 2180}%
\special{pa 4370 1550}%
\special{dt 0.027}%
\special{pa 5000 2240}%
\special{pa 4370 1610}%
\special{dt 0.027}%
\special{pa 5000 2300}%
\special{pa 4370 1670}%
\special{dt 0.027}%
\special{pa 5000 2360}%
\special{pa 4370 1730}%
\special{dt 0.027}%
\special{pa 5000 2420}%
\special{pa 4370 1790}%
\special{dt 0.027}%
\special{pa 5000 2480}%
\special{pa 4370 1850}%
\special{dt 0.027}%
\special{pa 5000 2540}%
\special{pa 4370 1910}%
\special{dt 0.027}%
\special{pa 5000 2600}%
\special{pa 4370 1970}%
\special{dt 0.027}%
\special{pa 5000 2660}%
\special{pa 4370 2030}%
\special{dt 0.027}%
\special{pa 5000 2720}%
\special{pa 4370 2090}%
\special{dt 0.027}%
\special{pa 5000 2780}%
\special{pa 4370 2150}%
\special{dt 0.027}%
\special{pa 4990 2830}%
\special{pa 4370 2210}%
\special{dt 0.027}%
\special{pa 4930 2830}%
\special{pa 4370 2270}%
\special{dt 0.027}%
\special{pa 4870 2830}%
\special{pa 4370 2330}%
\special{dt 0.027}%
\special{pa 4810 2830}%
\special{pa 4370 2390}%
\special{dt 0.027}%
\special{pa 4750 2830}%
\special{pa 4370 2450}%
\special{dt 0.027}%
\special{pa 4690 2830}%
\special{pa 4370 2510}%
\special{dt 0.027}%
\special{pa 4630 2830}%
\special{pa 4370 2570}%
\special{dt 0.027}%
\special{pa 4570 2830}%
\special{pa 4370 2630}%
\special{dt 0.027}%
\special{pa 4510 2830}%
\special{pa 4370 2690}%
\special{dt 0.027}%
\special{pa 4450 2830}%
\special{pa 4370 2750}%
\special{dt 0.027}%
\special{pa 5050 2830}%
\special{pa 5000 2780}%
\special{dt 0.027}%
\special{pa 5110 2830}%
\special{pa 5000 2720}%
\special{dt 0.027}%
% LINE 3 2 3 1
% 60 5170 2830 5000 2660 5230 2830 5000 2600 5290 2830 5000 2540 5350 2830 5000 2480 5410 2830 5000 2420 5470 2830 5000 2360 5530 2830 5000 2300 5590 2830 5000 2240 5650 2830 5000 2180 5710 2830 5000 2120 5770 2830 5000 2060 5830 2830 5000 2000 5890 2830 5000 1940 5950 2830 5000 1880 6010 2830 5000 1820 6070 2830 5000 1760 6130 2830 5000 1700 6190 2830 5000 1640 6250 2830 5000 1580 6310 2830 5000 1520 6370 2830 5000 1460 6410 2810 5000 1400 6410 2750 5000 1340 6410 2690 5000 1280 6410 2630 5000 1220 6410 2570 5050 1210 6410 2510 5110 1210 6410 2450 5170 1210 6410 2390 5230 1210 6410 2330 5290 1210
% 
\special{pn 4}%
\special{pa 5170 2830}%
\special{pa 5000 2660}%
\special{dt 0.027}%
\special{pa 5230 2830}%
\special{pa 5000 2600}%
\special{dt 0.027}%
\special{pa 5290 2830}%
\special{pa 5000 2540}%
\special{dt 0.027}%
\special{pa 5350 2830}%
\special{pa 5000 2480}%
\special{dt 0.027}%
\special{pa 5410 2830}%
\special{pa 5000 2420}%
\special{dt 0.027}%
\special{pa 5470 2830}%
\special{pa 5000 2360}%
\special{dt 0.027}%
\special{pa 5530 2830}%
\special{pa 5000 2300}%
\special{dt 0.027}%
\special{pa 5590 2830}%
\special{pa 5000 2240}%
\special{dt 0.027}%
\special{pa 5650 2830}%
\special{pa 5000 2180}%
\special{dt 0.027}%
\special{pa 5710 2830}%
\special{pa 5000 2120}%
\special{dt 0.027}%
\special{pa 5770 2830}%
\special{pa 5000 2060}%
\special{dt 0.027}%
\special{pa 5830 2830}%
\special{pa 5000 2000}%
\special{dt 0.027}%
\special{pa 5890 2830}%
\special{pa 5000 1940}%
\special{dt 0.027}%
\special{pa 5950 2830}%
\special{pa 5000 1880}%
\special{dt 0.027}%
\special{pa 6010 2830}%
\special{pa 5000 1820}%
\special{dt 0.027}%
\special{pa 6070 2830}%
\special{pa 5000 1760}%
\special{dt 0.027}%
\special{pa 6130 2830}%
\special{pa 5000 1700}%
\special{dt 0.027}%
\special{pa 6190 2830}%
\special{pa 5000 1640}%
\special{dt 0.027}%
\special{pa 6250 2830}%
\special{pa 5000 1580}%
\special{dt 0.027}%
\special{pa 6310 2830}%
\special{pa 5000 1520}%
\special{dt 0.027}%
\special{pa 6370 2830}%
\special{pa 5000 1460}%
\special{dt 0.027}%
\special{pa 6410 2810}%
\special{pa 5000 1400}%
\special{dt 0.027}%
\special{pa 6410 2750}%
\special{pa 5000 1340}%
\special{dt 0.027}%
\special{pa 6410 2690}%
\special{pa 5000 1280}%
\special{dt 0.027}%
\special{pa 6410 2630}%
\special{pa 5000 1220}%
\special{dt 0.027}%
\special{pa 6410 2570}%
\special{pa 5050 1210}%
\special{dt 0.027}%
\special{pa 6410 2510}%
\special{pa 5110 1210}%
\special{dt 0.027}%
\special{pa 6410 2450}%
\special{pa 5170 1210}%
\special{dt 0.027}%
\special{pa 6410 2390}%
\special{pa 5230 1210}%
\special{dt 0.027}%
\special{pa 6410 2330}%
\special{pa 5290 1210}%
\special{dt 0.027}%
% LINE 3 2 3 2
% 52 6410 2270 5350 1210 6410 2210 5410 1210 6410 2150 5470 1210 6410 2090 5530 1210 6410 2030 5590 1210 6410 1970 5650 1210 6410 1910 5710 1210 6410 1850 5770 1210 6410 1790 5830 1210 6410 1730 5890 1210 6410 1670 5950 1210 6410 1610 6010 1210 6410 1550 6070 1210 6410 1490 6130 1210 6410 1430 6190 1210 6410 1370 6250 1210 6410 1310 6310 1210 6410 1250 6370 1210 5000 1700 4510 1210 5000 1640 4570 1210 5000 1580 4630 1210 5000 1520 4690 1210 5000 1460 4750 1210 5000 1400 4810 1210 5000 1340 4870 1210 5000 1280 4930 1210
% 
\special{pn 4}%
\special{pa 6410 2270}%
\special{pa 5350 1210}%
\special{dt 0.027}%
\special{pa 6410 2210}%
\special{pa 5410 1210}%
\special{dt 0.027}%
\special{pa 6410 2150}%
\special{pa 5470 1210}%
\special{dt 0.027}%
\special{pa 6410 2090}%
\special{pa 5530 1210}%
\special{dt 0.027}%
\special{pa 6410 2030}%
\special{pa 5590 1210}%
\special{dt 0.027}%
\special{pa 6410 1970}%
\special{pa 5650 1210}%
\special{dt 0.027}%
\special{pa 6410 1910}%
\special{pa 5710 1210}%
\special{dt 0.027}%
\special{pa 6410 1850}%
\special{pa 5770 1210}%
\special{dt 0.027}%
\special{pa 6410 1790}%
\special{pa 5830 1210}%
\special{dt 0.027}%
\special{pa 6410 1730}%
\special{pa 5890 1210}%
\special{dt 0.027}%
\special{pa 6410 1670}%
\special{pa 5950 1210}%
\special{dt 0.027}%
\special{pa 6410 1610}%
\special{pa 6010 1210}%
\special{dt 0.027}%
\special{pa 6410 1550}%
\special{pa 6070 1210}%
\special{dt 0.027}%
\special{pa 6410 1490}%
\special{pa 6130 1210}%
\special{dt 0.027}%
\special{pa 6410 1430}%
\special{pa 6190 1210}%
\special{dt 0.027}%
\special{pa 6410 1370}%
\special{pa 6250 1210}%
\special{dt 0.027}%
\special{pa 6410 1310}%
\special{pa 6310 1210}%
\special{dt 0.027}%
\special{pa 6410 1250}%
\special{pa 6370 1210}%
\special{dt 0.027}%
\special{pa 5000 1700}%
\special{pa 4510 1210}%
\special{dt 0.027}%
\special{pa 5000 1640}%
\special{pa 4570 1210}%
\special{dt 0.027}%
\special{pa 5000 1580}%
\special{pa 4630 1210}%
\special{dt 0.027}%
\special{pa 5000 1520}%
\special{pa 4690 1210}%
\special{dt 0.027}%
\special{pa 5000 1460}%
\special{pa 4750 1210}%
\special{dt 0.027}%
\special{pa 5000 1400}%
\special{pa 4810 1210}%
\special{dt 0.027}%
\special{pa 5000 1340}%
\special{pa 4870 1210}%
\special{dt 0.027}%
\special{pa 5000 1280}%
\special{pa 4930 1210}%
\special{dt 0.027}%
% BOX 2 1 3 0
% 2 2580 1200 5030 3620
% 
\special{pn 8}%
\special{pa 2580 1200}%
\special{pa 5030 1200}%
\special{pa 5030 3620}%
\special{pa 2580 3620}%
\special{pa 2580 1200}%
\special{da 0.070}%
\end{picture}
\caption{A bad way of dividing $\Omega$}
\end{figure}
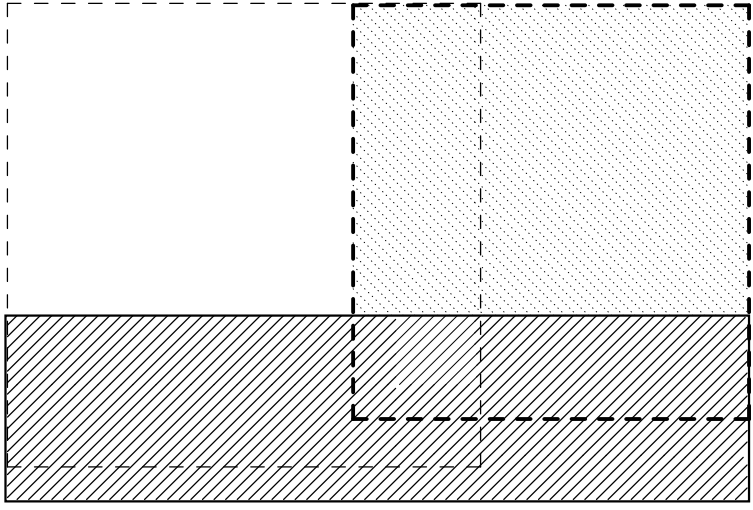
Figure 1 gives an example which satisfies our assumptions about the way we divide $\Omega$ into several subdomains. In Figure 2, since there is an overlapping area between the three subdomains, this way of dividing $\Omega$ does not satisfy our conditions.
\end{remark}
The Schwarz waveform relaxation algorithm solves $I$ equations in $I$ subdomains instead of solving directly the main problem $(\ref{2e1})$. The iterate $\#k$ in the $l$-th domain, denoted by $u_{l}^{k}$, is defined by
\begin{equation}\label{2e2}
\left \{
\begin{array}{ll}
{\partial_t u^k_l}-\sum_{i,j=1}^n\partial_j(a_{i,j}{\partial_{i}  u^k_l})+\sum_{i=1}^nb_{i}{\partial_i  u^k_l}\vspace{.1in}+c u^k_l=F(t,x,u^k_l),\mbox{ in }\Omega_l\times(0,\infty),\vspace{.1in}\\
 \mathfrak{B}_{l,l'}u^k_l=\mathfrak{B}_{l,l'}u^{k-1}_{l'}, \mbox{ on }\Gamma_{l,l'}\times(0,\infty),\forall l'\in J_l,\end{array}\right. 
\end{equation}
where the transmission operator $\mathfrak{B}_{l,l'}$ is either of the Dirichlet type or of the Robin type.
\\ Each iterate inherits the boundary conditions and the initial values of $u$
\[
u_l^k= g \mbox{ on }(\partial\Omega_l\cap\partial\Omega)\times(0,\infty),\quad
u_l^k(.,0)= g(.,0) \mbox{ in } \Omega_l.
\]
A bounded initial guess $u^0$ in $C^{\infty}(\overline{\Omega\times(0,\infty)})$ is provided, \textit{i.e.} at step $1$ Equations $(\ref{2e2})$ are solved
\begin{eqnarray*}
\mathfrak{B}_{l,l'}u^1_l=u^{0} &\mbox{ on }\Gamma_{l,l'}\times(0,\infty),\forall l'\in J_l.
\end{eqnarray*}
We assume also the compatibility condition on $u^{0}$
\begin{eqnarray*}
\mathfrak{B}_{l,l'}g(.,0)=u^{0}(.,0) &\mbox{ on }\Gamma_{l,l'},\forall l'\in J_l.
\end{eqnarray*}
\\ Denote by $e_l^k$ the difference between $u_l^k$ and $u$, and subtract Equation $(\ref{2e2})$ with the main equation $(\ref{2e1})$ to obtain the following equations on $e_l^k$
\begin{equation}\label{2e3}
\left \{
\begin{array}{ll}
{\partial_t e^k_l}-\sum_{i,j=1}^n\partial_j(a_{i,j}{\partial_{i}  e^k_l})+\sum_{i=1}^nb_{i}{\partial_t  e^k_l}+c e^k_l=F(t,x,u^k_l)-F(t,x,u) \mbox{ in }\Omega_l\times(0,\infty),\vspace{.1in}\\
 \mathfrak{B}_{l,l'}e^k_l=\mathfrak{B}_{l,l'}e^{k-1}_{l'} \mbox{ on }\Gamma_{l,l'}\times(0,\infty),\forall l'\in J_l.\end{array}\right. 
\end{equation}
Moreover, 
\[
e_l^k(.,.)= 0 \mbox{ on }(\partial\Omega_l\cap\partial\Omega)\times(0,\infty),\quad
e_l^k(.,0)= 0 \mbox{ in } \Omega_l.
\]
\subsection{Classical Schwarz Methods}
Consider the classical Schwarz waveform relaxation algorithm with Dirichlet transmission conditions $\mathfrak{B}_{l,l'}=Id$. By induction, each subproblem $(\ref{2e2})$ in each iteration has a unique solution  in $C^{2,1}(\overline{\Omega\times(0,\infty)})$ then in $L^2(0,\infty,H^1(\Omega))\cap L^{\infty}(\Omega\times(0,\infty))$ also. 
\\ Consider $(\ref{2e3})$ and let $g$, $f$ be bounded and strictly positive functions in $C^{\infty}(\mathbb{R}^n,\mathbb{R})$ and $C^{\infty}(\mathbb{R},\mathbb{R})$. Define 
\begin{eqnarray*}
\Phi_l^k(x,t):=(e_l^k)^2g(x)f(t).
\end{eqnarray*}
Since $e_l^k$ belongs to $L^2(0,\infty,H^1(\Omega))\cup L^{\infty}(\Omega\times(0,\infty))$, $\Phi_l^k$ belongs to $L^2(0,\infty,H^1(\Omega))$.
\\ Let $c_i$ be $b_i+\sum_{j=1}^n2a_{i,j}{\partial_j g}g^{-1}$, then $c_i\in L^{\infty}(\Omega_l\times(0,\infty))$, and define the following operator
\begin{eqnarray}\label{2e11}
\mathfrak{L}_{lD}(\Phi)&=&{\partial_t\Phi}-\sum_{i,j=1}^n\partial_j(a_{i,j}{\partial_{i,j}\Phi})+\sum_{i=1}^{n}c_i(x,t){\partial_i\Phi}.
\end{eqnarray}
\begin{lemma}\label{linear1}
 In each subdomain $\Omega_l$, for each iterate $k$ 
$$\mathfrak{L}_{lD}(\Phi_l^k)\leq0,$$
in the distributional sense, {\it i.e.} for all $\varphi$ in $H_0^1(\Omega)$ and $\varphi\geq0$ a.e. on $\Omega$,  $$\int_{\Omega_l}\mathfrak{L}_{lD}(\Phi_l^k)\varphi\leq0 \mbox{ a.e. in }(0,\infty).$$
\end{lemma}
\begin{proof} Define the operator
\begin{eqnarray*}
\mathfrak{L}_{lD0}:={\partial_t}-\sum_{i,j=1}^n\partial_j(a_{i}{\partial_{i}}).
\end{eqnarray*}
%\\\h  Derive the function $\Phi_l^k$
%\begin{eqnarray}\label{2e4}
%{\partial_t\Phi_l^k}& = &2{\partial_t e_l^k}e_l^kgf+(e_l^k)^2gf', \\\nonumber
%{\partial_i\Phi_l^k}& = &2{\partial_i e_l^k}e_l^kgf+(e_l^k)^2{\partial_i g}f, \forall i\in\{1,n\},\\\nonumber
%{\partial_{i,j}\Phi_l^k}& = &2{\partial_{i,j} e_l^k}e_l^kgf+2{\partial_i e_l^k}{\partial_j %e_l^k}gf+2{\partial_i e_l^k}e_l^k{\partial_i g}f\\\nonumber
%& &+2e_l^k{\partial_j e_l^k}{\partial_i g}f+(e_l^k)^2{\partial_{i,j} g}f,\forall i,j\in\{1,n\}.
%\end{eqnarray}
A lengthy but easy computation then implies that 
\begin{eqnarray}\label{2e5}
\mathfrak{L}_{lD0}(\Phi_l^k)& =&2\left({\partial_t e_l^k}-\sum_{i,j=1}^n\partial_j(a_{i,j}{\partial_{i,j}e_l^k})\right)e_l^kgf-\sum_{i,j=1}^n2a_{i,j}{\partial_{i} e_l^k}{\partial_j e_l^k}gf-\\\nonumber
& &-\sum_{i,j=1}^n4a_{i,j}{\partial_i e_l^k}e_l^k{\partial_j g}f+(e_l^k)^2\left(-\sum_{i,j=1}^n(a_{i,j}{\partial_{i,j} g}+\partial_i a_{i,j}\partial_j g)f+gf'\right).
\end{eqnarray}
Thanks to Equation $(\ref{2e3})$, and the lipschitzian property of $F$, the first term in $(\ref{2e5})$ can be estimated in the distributional sense
\begin{eqnarray}\label{2e6}
& &2\left({\partial_t e_l^k}-\sum_{i,j=1}^n\partial_j(a_{i,j}{\partial_{i}e_l^k})\right)e_l^kgf\\\nonumber
&=&2gfe_l^k\left(F(t,x,u_l^k)-F(t,x,u)-\sum_{i=1}^nb_i{\partial_i e_l^k}-ce_l^k\right)\\\nonumber
&\leq&2gf(e_l^k)^2(C+||c||_{\infty})-\sum_{i=1}^n2b_i{\partial_i e_l^k}gfe_l^k\\\nonumber
&\leq&2gf(e_l^k)^2(C+||c||_{\infty})-\sum_{i=1}^nb_i\left({\partial_i\Phi_l^k}-(e_l^k)^2{\partial_i g}f\right)\\\nonumber
&\leq&-\sum_{i=1}^nb_i{\partial_i\Phi_l^k}+(e_l^k)^2\left(2gf(C+||c||_{\infty})+\sum_{i=1}^{n}b_i{\partial_i g}f\right).
\end{eqnarray}
Since $\mathfrak{L}_{lD0}$ is elliptic, the second term in $(\ref{2e5})$ is negative.
\\ Moreover, the third term in $(\ref{2e5})$ can be transformed into
\begin{eqnarray}\label{2e7}
-\sum_{i,j=1}^n4a_{i,j}{\partial_i e_l^k}e_l^k{\partial_j g}f
%&=&-\sum_{i,j=1}^n2a_{i,j}{\partial_j g}g^{-1}\left({\partial_i\Phi_l^k}-(e_l^k)^2{\partial_i g}f\right)\\\nonumber
=-\sum_{i,j=1}^n2a_{i,j}{\partial_j g}g^{-1}{\partial_i\Phi_l^k}+\sum_{i,j=1}^n2a_{i,j}{\partial_j g}{\partial_i g}g^{-1}f(e_l^k)^2.\nonumber
\end{eqnarray}
Combine $(\ref{2e5})$, $(\ref{2e6})$ and $(\ref{2e7})$ to get
\begin{eqnarray}\label{2e9}
\mathfrak{L}_{lD0}(\Phi_l^k)+\sum_{i=1}^{n}\left(b_i+\sum_{j=1}^n2a_{i,j}{\partial_j g}g^{-1}\right){\partial_i\Phi_l^k}\leq(e_l^k)^2\mathfrak{M},
\end{eqnarray}
where 
\begin{eqnarray}\label{2e10}
\mathfrak{M}%&=&-\sum_{i,j=1}^na_{i,j}{\partial_{i,j} g}f+gf'+2gf(C+||c||_{\infty})+\\\nonumber
%& &+\sum_{i=1}^{n}b_i{\partial_i g}f+\sum_{i,j=1}^n2a_{i,j}{\partial_j g}{\partial_i g}g^{-1}f\\\nonumber
=\left[\sum_{i,j=1}^n(-a_{i,j}\frac{{\partial_{i,j} g}}{g}-\partial_ja_{i,j}\frac{\partial_i g}{g}+2a_{i,j}\frac{{\partial_j g}}{g}\frac{{\partial_i g}}{g})+\frac{f'}{f}+2(C+||c||_{\infty})+\sum_{i=1}^{n}b_i\frac{{\partial_i g}}{g}\right]fg.
\end{eqnarray}
\\\h Notice that $\mathfrak{M}$ has the form of $fg(\frac{f'}{f}+G(g))$. Choosing $f$ such that $-\frac{f'}{f}$ large enough (for example, $f=\exp(-\alpha t)$, where $\alpha$ is a large positive constant), since the other terms are bounded in the bounded domain $\Omega_l\times(0,\infty)$, then $\mathfrak{M}\leq 0$. The nonlinear equation $(\ref{2e3})$ has been transformed into the following linearized inequation  of $\Phi_l^k$ 
\begin{eqnarray}\label{2e11}
{\partial_t\Phi_l^k}-\sum_{i,j=1}^n\partial_j(a_{i,j}{\partial_{i}\Phi_l^k})+\sum_{i=1}^{n}c_i{\partial_i\Phi_l^k}\leq0.
\end{eqnarray}
\end{proof}
\begin{theorem}\label{2t1} Consider the Schwarz algorithm with Dirichlet transmission condition, suppose that $f(t)$ is a strictly positive and continuous function satisfying 
$-\min_{t\in(0,\infty)}\frac{f'(t)}{f(t)}$ is sufficiently large, we get the geometrical convergence
\begin{eqnarray*}
\mathop{\lim}_{k\to\infty}\max_{l\in\{1,\dots, I\}}||(u_l^k-u)^2f(t)||_{L^{\infty}(\Omega_l\times(0,\infty))}=0.
\end{eqnarray*}
\end{theorem}
\begin{remark}
In the above theorem, if $f$ is chosen to be $\exp(-\alpha t)$, then when $\alpha$ is large enough, $$-\min_{t\in(0,\infty)}\frac{f'(t)}{f(t)}$$ is large enough. In this case, the limit
\begin{eqnarray*}
\mathop{\lim}_{k\to\infty}\max_{l\in\{1,\dots, I\}}||(u_l^k-u)^2f(t)||_{L^{\infty}(\Omega_l\times(0,\infty))}=0
\end{eqnarray*}
implies the almost everywhere convergence of the sequence $\{u_l^k\}$ to $u$ on $\Omega_l\times(0,\infty)$.
\end{remark}
\begin{remark}
In the proof, if $a_{ij}$, $b_i$, $c$ depend both on $t$ and $x$, the convergence result in the theorem remains true.
\end{remark}
\begin{proof} The proof is divided into two steps.
\\{\bf Step 1:} Construct estimates of the errors $e_l^k$ from Inequation $(\ref{2e11})$.
\\\h Define 
\begin{eqnarray}\label{2e12}
M=\mathop{esssup}_{\partial\Omega_l\times[0,\infty)\cup\Omega_l\times\{0\}}\Phi_l^k(x,t),
\end{eqnarray}
we will prove that the maximum principle holds, i.e. $M\geq \Phi_l^k$ a.e. on $\Omega_l\times(0,\infty)$. Define the function $$w=(\Phi_l^k-M)_+=\max\{\Phi_l^k-M,0\}.$$ Since $w\in L^2(0,\infty,H_0^{1}(\Omega_l)),$ then
\begin{eqnarray}\label{2e13}
{\partial_t w}-\sum_{i,j=1}^n\partial_j(a_{i,j}{\partial_{i} w})+\sum_{i=1}^{n}c_i{\partial_i w}\leq0.
\end{eqnarray}
To prove that $M\geq \Phi_l^k$ a.e. on $\Omega_l\times(0,\infty)$, it suffices to prove that $w=0$ a.e. on $\Omega_l\times(0,\infty)$.
\\\h For $0<p\leq\infty$ and $0<q\leq\infty$, and for $0\leq\tau_1<\tau_2\leq\infty$, define $||h||_{\tau_1,\tau_2,p,q}=||h||_{L^q(\tau_1,\tau_2,L^p(\Omega_l))}$, for $h\in L^q(\tau_1,\tau_2,L^p(\Omega_l))$. If $\tau_1=0$, denote $||h||_{0,\tau_2,p,q}$ by $||h||_{\tau_2,p,q}$.
\\\h Let $\chi(\tau_1,\tau_2)$ be the characteristic function of the open interval $(\tau_1,\tau_2)$, where $0<\tau_1<\tau_2\leq \infty$ and set $\varphi=w\chi$. Since $w\in L^2(0,T,H_0^{1}(\Omega_l))\cap L^\infty(\Omega_l\times(0,\infty))$, it is evident that $\varphi\in  L^2(0,\infty,H_0^{1}(\Omega_l))\cap L^\infty(\Omega_l\times(0,\infty))$.
\\\h Use $\varphi$ as a test function for $(\ref{2e13})$ 
\begin{eqnarray}\label{2e14}
\int_{\tau_1}^{\tau_2}\int_{\Omega_l}{\partial_t w}w dxdt+\int_{\tau_1}^{\tau_2}\int_{\Omega_l}\sum_{i,j=1}^na_{i,j}{\partial_i w}{\partial_j w}dxdt+\int_{\tau_1}^{\tau_2}\int_{\Omega_l}\sum_{i=1}^{n}c_i{\partial_i w}w dxdt\leq0.
\end{eqnarray}
Equation $(\ref{2e14})$ and Conditions $(A_1)$ and $(A_2)$ imply  that
\begin{eqnarray}\label{2e15}
& &\left.\int_{\Omega_l} \frac{w^{2}}{2}dx\right|_{t=\tau_1}^{t=\tau_2}+\int_{\tau_1}^{\tau_2}\int_{\Omega_l}\lambda|\nabla w|^2 dxdt\leq\int_{\tau_1}^{\tau_2}\int_{\Omega_l}M_1|\nabla w|w dxdt,
\end{eqnarray}
where $M_1$ is a positive constant. By the Cauchy inequality, the right hand side of $(\ref{2e15})$ is bounded by
\begin{eqnarray}\label{2e16}
\int_{\tau_1}^{\tau_2}\int_{\Omega_l}\frac{M_1\epsilon}{2}|\nabla w|^2 dxdt+\int_{\tau_1}^{\tau_2}\int_{\Omega_l}\frac{M_1}{2\epsilon}w^{2} dxdt,
\end{eqnarray}
where $\epsilon$ is a small positive constant.
\\ For $\epsilon$ to be $\frac{\lambda}{M_1}$, Equality $(\ref{2e16})$ implies
\begin{eqnarray}\label{2e17}
\left.\int_{\Omega_l}\frac{w^{2}}{2}dx\right|_{t=\tau_1}^{t=\tau_2}+\int_{\tau_1}^{\tau_2}\int_{\Omega_l}\frac{\lambda}{2} |\nabla w|^2dxdt\leq\int_{\tau_1}^{\tau_2}\int_{\Omega_l}M_2w^{2} dxdt,
\end{eqnarray}
with $M_2=\frac{M_1^2}{2}$.
\\\h Denote $X(t)=\int_{\Omega_l}w^{2}(x,t)dx$, and let $\tau_2$ be $t$ in the interval $\mathfrak{I}=(\tau_1,\tau_1+\delta)$, the previous estimate infers
\begin{eqnarray}\label{2e18}
X(t)+\lambda||\nabla w||^2_{\tau_1,t,2,2}\leq M_2\delta\{\sup_{\mathfrak{I}}X(t)\}+X(\tau_1).
\end{eqnarray}
Choosing $\delta$ such that $M_{2}\delta=\frac{1}{2}$, the fact that $\sup_{\mathfrak{I}}X(t)\leq 2X(\tau_1)$ then follows. Since the inequality is true on any time interval with the length of $\delta$, and $X(0)=0$, then $X(t)=0$ for a.e. $t$ in $(0,\infty)$. Hence $w=0$ for a.e. $t$ in $(0,\infty)$.
\\\h We have just proved that 
\begin{eqnarray}\label{2e19}
\left(e_l^k(x,t)\right)^2 g(x)f(t)\leq \max_{l'\in J_l}\left(\mathop{esssup}_{\Gamma_{l,l'}\times(0,\infty)}\left(e_{l}^k(x,t)\right)^2g(x)f(t)\right),
\end{eqnarray}
for all $l$ in $I$, for a.e. $(x,t)$ in $\Omega_l\times(0,\infty)$; for any strictly positive functions $g$, $f$  in $C^{\infty}(\mathbb{R}^n,\mathbb{R})$ and $C^{\infty}(\mathbb{R},\mathbb{R})$. 
{\\\bf Step 2:} The convergence of the algorithm.
\\\h Denote
\begin{eqnarray}\label{2e20}
E^k=\max_{l\in I}\left(\mathop{esssup}_{(x,t)\in(\Omega_l\times(0,\infty))}\left(e_l^k\right)^2 f(t)\right).
\end{eqnarray}
\h From $(\ref{2e19})$ comes that for every $l'$ in $J_l$ and for a.e. $(x,t)$ in $ {\Gamma_{l,l'}\times(0,\infty)}$
\begin{eqnarray}\label{2e21}
\left(e_{l}^k(x,t)\right)^2g(x)f(t)\leq \max_{l''\in J_{l'}}\left(\mathop{esssup}_{\Gamma_{l',l''}\times(0,\infty)}\left(e_{l'}^{k-1}(x,t)\right)^2g(x)f(t)\right),
\end{eqnarray}
that implies
\begin{eqnarray}\label{2e22}
\left(e_{l}^k(x,t)\right)^2f(t)\leq \frac{1}{g(x)}\max_{l''\in J_{l'}}\left(\mathop{esssup}_{\Gamma_{l',l''}\times(0,\infty)}\left(e_{l''}^{k-2}(x,t)\right)^2g(x)f(t)\right).
\end{eqnarray}
Since $\Gamma_{l,l'}$ lies inside $\Omega_{l'}$, choose $g$ such that there exists a constant $M_{3,l}$ satisfying
$$1>M_{3,l}>\frac{g(\zeta')}{g(\zeta)},~~~\forall \zeta\in \Gamma_{l,l'} \mbox{ and }\forall\zeta'\in \cup_{l''\in J_{l'}}\Gamma_{l',l''},$$ and this implies  that   for all $l'$ in $J_l$, for a.e. $(x,t)$ in ${\Gamma_{l,l'}\times(0,\infty)}$
\begin{eqnarray}\label{2e23}
\left(e_{l}^k(x,t)\right)^2f(t)\leq M_{3,l}\max_{l''\in J_{l'}}\left(\mathop{esssup}_{\Gamma_{l',l''}\times(0,\infty)}\left(e_{l''}^{k-2}(x,t)\right)^2f(t)\right)\leq M_{3,l}E^{k-2}.
\end{eqnarray}
\\ Choose $g$ to be the function $1$, $(\ref{2e19})$ yields
\begin{eqnarray}\label{2e24}
\left(e_l^k(x,t)\right)^2 f(t)\leq \max_{l'\in J_l}\left(\mathop{esssup}_{\Gamma_{l,l'}\times(0,\infty)}\left(e_{l}^k(x,t)\right)^2f(t)\right),~~~\forall l'\in J_l, \mbox{ a.e. on }\Omega_l.
\end{eqnarray}
\\\h The estimates $(\ref{2e23})$ and $(\ref{2e24})$ imply the existence of a constant $M_4$ smaller than 1 and satisfy
\begin{eqnarray}\label{2e25}
E^k\leq M_4E^{k-2},
\end{eqnarray}
which shows that the errors converge geometrically $$\lim_{k\to\infty}E^k=0.$$ The theorem is proved.
\end{proof}
\subsection{Optimized Schwarz Methods}
The optimized Schwarz waveform relaxation algorithms  are defined by replacing the Dirichlet by Robin transmission operators
$$\mathfrak{B}_{l,l'}v=\sum_{i,j=1}^n a_{i,j}{\partial_i v}n_{l,l',j}+p_{l,l'}v,$$ 
where $n_{l,l',j}$ is the $j$-th component of the outward unit normal vector of $\Gamma_{l,l'}$; $p_{l,l'}$ is positive and belongs to $L^{\infty}(\Gamma_{l,l'})$.
By induction, each subproblem $(\ref{2e2})$ in each iteration has a unique solution in $L^{2}(0,\infty,H^1(\Omega))$ and the algorithm is well-posed. 
\\ Let $f$ be a function in $L^2(0,\infty)$, define
\begin{equation*}
\int_0^{\infty}f(x)\exp(-yx)dx.
\end{equation*}
Now, define for a fixed positive number $\alpha$
\begin{equation*}
|f|_{\alpha}=\sup_{\alpha'>\alpha}\left[\int_{\alpha'}^{\alpha'+1}\left(\int_0^{\infty}f(x)\exp(-y x)dx\right)^2dy\right]^{\frac{1}{2}},
\end{equation*}
and
$$\mathbb{L}_{\alpha}^{2}(0,\infty)=\{f~:~f\in L^2(0,\infty), |f|_{\alpha}<\infty\}.$$
Then $(\mathbb{L}_{\alpha}^{2}(0,\infty),|.|_\alpha)$ is a normed subspace of $L^2(0,\infty)$.
%\begin{lemma}
%$(\mathbb{L}_{\alpha}^{2}(0,\infty),|.|_\alpha)$ is a new normed subspace of $L^2(0,\infty)$.
%\end{lemma}
%\begin{proof}
%Thus, $|f|_{\alpha}=0$ if and only if $f=0$ (see, for example the book \cite{Widder:1941:TLT}). Moreover, for any $\gamma$ in $\mathbb{R}$, $|\gamma f|_{\alpha}=|\gamma||f|_{\alpha}$. Consider $f$, $g$ in $L^2(0,\infty)$ such that $|f|_{\alpha}$, $|g|_{\alpha}$ $<$ $\infty$, by Holder's inequality, for $\alpha'$ being greater than $\alpha$
%\begin{equation*}
%\left[\int_{\alpha'}^{\alpha'+1}\left(\int_0^{\infty}f(x)\exp(-y x)dx\right)^2dy\right]^{\frac{1}{2}}\geq \int_{\alpha'}^{\alpha'+1}\left|\int_0^{\infty}f(x)\exp(-y x)dx\right|dy,
%\end{equation*}
%and
%\begin{equation*}
%\left[\int_{\alpha'}^{\alpha'+1}\left(\int_0^{\infty}g(x)\exp(-y x)dx\right)^2dy\right]^{\frac{1}{2}}\geq %\int_{\alpha'}^{\alpha'+1}\left|\int_0^{\infty}g(x)\exp(-y x)dx\right|dy.
%\end{equation*}
%These estimates lead to the following 
%\begin{eqnarray*}
%& &\left[\int_{\alpha'}^{\alpha'+1}\left(\int_0^{\infty}(f+g)\exp(-y x)dx\right)^2dy\right]^{\frac{1}{2}}\\
%& \leq & \left[\int_{\alpha'}^{\alpha'+1}\left(\int_0^{\infty}f\exp(-y x)dx\right)^2dy\right]^{\frac{1}{2}}+\left[\int_{\alpha'}^{\alpha'+1}\left(\int_0^{\infty}g\exp(-y x)dx\right)^2dy\right]^{\frac{1}{2}}.
%\end{eqnarray*}
%Hence $$|f+g|_{\alpha}\leq |f|_\alpha+|g|_\alpha.$$
%\end{proof}
\\\h Consider Equation $(\ref{2e3})$, let $g_l$ be a function bounded and greater than $1$  in $C^{\infty}(\mathbb{R}^n,\mathbb{R})$, $\alpha$ be a positive constant, and define 
\begin{eqnarray*}
\Phi_l^k(x):=\left(\int_0^{\infty}e_l^k\exp(-\alpha t)dt\right)g_l(x),
\end{eqnarray*}
then $\Phi_l^k(x)$ belongs to $H^1(\Omega_l)$.
\\\h Let $B^l_{i}$ and $C^l$ be functions in $L^{\infty}(\mathbb{R}^n)$ defined in the following ways:
$$B^l_{i}:=b_i+\sum_{j=1}^n\left(a_{i,j}\frac{\partial_j g_l}{g_l}\right),$$
$$C^l=\left[\frac{\alpha}{2}+\sum_{i,j=1}^n\left(-a_{i,j}\frac{2{\partial_i g_l}{\partial_j g_l}}{(g_l)^2}-\partial_ja_{i,j}\frac{\partial_i g}{g}+a_{i,j}\frac{{\partial_{i,j} g_l}}{g_l}\right)-\sum_{i=1}^nb_{i}\frac{{\partial_i g_l}}{g_l}\right].$$
Define
\begin{eqnarray}\label{2e32}
\mathfrak{L}_{lR}(\Phi_l^k)&=&-\sum_{i,j=1}^n\partial_j(a_{i,j}{\partial_{i} \Phi_l^k})+\sum_{i=1}^nB^l_{i}{\partial_i \Phi_l^k}+C^l\Phi_l^k\\\nonumber
& &+\left\{\int_0^{\infty}\left[\left(\frac{\alpha}{2}+c\right)e_l^k-F(u_l^k)+F(u)\right]\exp(-\alpha t)dt\right\}g_l.
\end{eqnarray}
 It is possible to suppose $\alpha$ to be large such that $C^l$ belongs to $(\frac{\alpha}{4},\alpha)$.
\begin{lemma}\label{linear2} Choose $g_l$, $g_{l'}$ such that $\nabla g_l=\nabla g_{l'}=0$ on $\Gamma_{l,l'}$ and $\frac{g_{l'}}{g_l}>1$ on $\Gamma_{l,l'}$, for all $l'$ in $J_l$.
$\Phi_l^k$ is then a solution of the following equation
\begin{equation}\label{2e34}\
\left \{
\begin{array}{ll}
\mathfrak{L}_{lR}(\Phi_l^k)=0,&\mbox{ in }\Omega_l\times(0,\infty),\vspace{.1in}\\
 \beta_l\mathfrak{B}_{l,l'}(\Phi_l^k)=\mathfrak{B}_{l,l'}(\Phi_{l'}^{k-1}) &\mbox{ on }\Gamma_{l,l'}\times(0,\infty),\forall l'\in J_l.\end{array}\right. 
\end{equation}
where $\beta_l=\frac{g_{l'}}{g_l}$ on $\Gamma_{l,l'}$, for all $l'$ in $J_l$.
\end{lemma}
\begin{proof}
%We have the following computation on the function, for all $i,j$ in $\{1,\dots,n\}$
%\begin{eqnarray}\label{2e26}
%& &\alpha\Phi_l^k=\\\nonumber
%&=&\left(\int_0^{\infty}{\partial_t e_l^k}\exp(-\alpha t)dt\right)g_l,\\\nonumber
%& &{\partial_i\Phi_l^k }=\\\nonumber
%&=&\left(\int_0^{\infty}{\partial_i e_l^k}\exp(-\alpha t)dt\right)g_l+\left(\int_0^{\infty}e_l^k\exp(-\alpha t)dt\right){\partial_i g_l},\\\nonumber
%& &{\partial_{i,j}\Phi_l^k }=\\\nonumber
%&=&\left(\int_0^{\infty}{\partial_j e_l^k}\exp(-\alpha t)dt\right){\partial_i g_l}+\left(\int_0^{\infty}%{\partial_i e_l^k}\exp(-\alpha t)dt\right){\partial_j g_l}\\\nonumber
%& &+\left(\int_0^{\infty}{\partial_{i,j} e_l^k}\exp(-\alpha t)dt\right)g_l+\left(\int_0^{\infty}e_l^k\exp(-\alpha t)dt\right){\partial_{i,j} g_l}.\\\nonumber
%\end{eqnarray}
A complicated but easy computation leads to 
\begin{eqnarray}\label{2e27}
& &-\sum_{i,j=1}^n\partial_j(a_{i,j}{\partial_{i}\Phi_l^k})+\alpha \Phi_l^k\\\nonumber
&=&\left[\int_0^{\infty}\left({\partial_t e_l^k}-\sum_{i,j=1}^n\partial_j(a_{i,j}{\partial_{i} e_l^k})\right)\exp(-\alpha t)dt\right]g_l\\\nonumber
& &-\left(\int_0^{\infty}e_l^k\exp(-\alpha t)dt\right)\left(\sum_{i,j=1}^na_{i,j}({\partial_{i,j} g_l}+\partial_ja_{i,j}\partial_i g_l)\right)\\\nonumber
& &-\sum_{i,j=1}^na_{i,j}\left[\int_0^{\infty}e_l^k\left({\partial_i e_l^k}{\partial_j g_l}+{\partial_j e_l^k}{\partial_i g_l}\right)\exp(-\alpha t)dt\right].
\end{eqnarray}
That implies
\begin{eqnarray}\label{2e28}
& &-\sum_{i,j=1}^n\partial_j(a_{i,j}{\partial_{i}\Phi_l^k})+\alpha \Phi_l^k\\\nonumber
& &+\left(\int_0^{\infty}e_l^k\exp(-\alpha t)dt\right)\left(\sum_{i,j=1}^n(a_{i,j}{\partial_{i,j} g_l}+\partial_ja_{i,j}\partial_i g_l)\right)\\\nonumber
&=&\left[\int_0^{\infty}\left(-\sum_{i=1}^nb_{i}{\partial_i e_l^k}-ce_l^k+F(u_l^k)-F(u)\right)\exp(-\alpha t)dt\right]g_l\\\nonumber
& &-\sum_{i,j=1}^na_{i,j}\left[\int_0^{\infty}e_l^k\left({\partial_i e_l^k}{\partial_j g_l}+{\partial_j e_l^k}{\partial_i g_l}\right)\exp(-\alpha t)dt\right].
\end{eqnarray}
%\\ Thanks to $(\ref{2e26})$, the first term on the right hand side of $(\ref{2e28})$ is equal to
%\begin{eqnarray}\label{2e29}
%& &-\sum_{i=1}^nb_{i}\left({\partial_i \Phi_l^k}-\int_0^{\infty}\rho(e_l^k)\exp(-\alpha t)dt{\partial_i g_l}\right)-\\\nonumber
%& &-\left[\int_0^{\infty}(ce_l^k-F(u_l^k)+F(u))\exp(-\alpha t)dt\right]g_l\\\nonumber
%& = &-\sum_{i=1}^nb_{i}{\partial_i \Phi_l^k}+\sum_{i=1}^nb_{i}\Phi_l^k\frac{{\partial_i g_l}}{g_l}-\left[\int_0^{\infty}(ce_l^k-F(u_l^k)+F(u))\exp(-\alpha t)dt\right]g_l.
%\end{eqnarray}
%Now, consider the second term on the right hand side of $(\ref{2e28})$, which is in fact
%\begin{eqnarray}\label{2e30}
%& &-\sum_{i,j=1}^na_{i,j}\left[\left({\partial_i \Phi_l^k}-\int_0^{\infty}e_l^k\exp(-\alpha t){\partial_i g_l}dt\right)\frac{{\partial_j g_l}}{g_l}+\right.\\\nonumber
%&&+\left.\left({\partial_j \Phi_l^k}-\int_0^{\infty}e_l^k\exp(-\alpha t){\partial_j g_l}dt\right)\frac{{\partial_i g_l}}{g_l}\right]\\\nonumber
%&=&-\sum_{i,j=1}^na_{i,j}\left({\partial_i \Phi_l^k}\frac{{\partial_j g_l}}{g_l}+{\partial_j \Phi_l^k}\frac{{\partial_i g_l}}{g_l}\right)+\sum_{i,j=1}^na_{i,j}\int_0^{\infty}e_l^k\exp(-\alpha t)dt\frac{2{\partial_i g_l}{\partial_j g_l}}{g_l}\\\nonumber
%&=&-\sum_{i,j=1}^na_{i,j}\left({\partial_i \Phi_l^k}\frac{{\partial_j g_l}}{g_l}+{\partial_j \Phi_l^k}\frac{{\partial_i g_l}}{g_l}\right)+\Phi_l^k\left(\sum_{i,j=1}^na_{i,j}\frac{2{\partial_i g_l}{\partial_j g_l}}{(g_l)^2}\right).
%\end{eqnarray}
%Combining $(\ref{2e28})$, $(\ref{2e29})$ and $(\ref{2e30})$, 
Therefore
\begin{eqnarray}\label{2e31}\nonumber
0&=&-\sum_{i,j=1}^n\partial_j(a_{i,j}{\partial_{i} \Phi_l^k})+\sum_{i=1}^nb_{i}{\partial_i \Phi_l^k}+\sum_{i,j=1}^na_{i,j}\left({\partial_i \Phi_l^k}\frac{{\partial_j g_l}}{g_l}+{\partial_j \Phi_l^k}\frac{{\partial_i g_l}}{g_l}\right)\\\nonumber
& &+\left(\frac{\alpha}{2}+\sum_{i,j=1}^na_{i,j}(-\frac{2{\partial_i g_l}{\partial_j g_l}}{(g_l)^2}+a_{i,j}\frac{{\partial_{i,j} g_l}}{g_l}-\partial_ja_{i,j}\frac{\partial_i g_l}{g_l})-\sum_{i=1}^nb_{i}\frac{{\partial_i g_l}}{g_l}\right)\Phi_l^k+\\\nonumber
& &+\left\{\int_0^{\infty}\left[\left(\frac{\alpha}{2}+c\right)e_l^k-F(u_l^k)+F(u)\right]\exp(-\alpha t)dt\right\}g_l.
\end{eqnarray}
\\\h Now, consider Robin transmission conditions on the boundary $\Gamma_{l,l'}$ and notice that  $\nabla g_l=\nabla g_{l'}=0$ on $\Gamma_{l,l'}$, the transmission conditions can be reported on $\Phi_l^k$
\begin{eqnarray}\label{2e33}
\mathfrak{B}_{l,l'}(\Phi_l^k)&=&\sum_{i,j=1}^na_{i,j}{\partial_i \Phi_l^k}n_{l,l',j}+p_{l,l'}\Phi_l^k\\\nonumber
&=&\int_0^{\infty}\left(\sum_{i,j=1}^na_{i,j}n_{l,l',j}{\partial_i e_l^k}+p_{l,l'}e_l^k\right)\exp(-\alpha t)g_l+\sum_{i,j=1}^na_{i,j}n_{l,l',j}{\partial_i g_l}e_l^k\\\nonumber
&=&\int_0^{\infty}\left(\sum_{i,j=1}^na_{i,j}n_{l,l',j}{\partial_i e_{l'}^{k-1}}+p_{l,l'}e_{l'}^{k-1}\right)\exp(-\alpha t)g_l\\\nonumber
&=&\left(\sum_{i,j=1}^na_{i,j}{\partial_i \Phi_{l'}^{k-1}}n_{l,l',j}+p_{l,l'}\Phi_{l'}^{k-1}\right)\frac{g_{l}}{g_{l'}}=\mathfrak{B}_{l,l'}(\Phi_{l'}^{k-1})\frac{g_{l}}{g_{l'}}.
\end{eqnarray}
\end{proof}
\begin{theorem}\label{2t2} Consider Schwarz algorithms with Robin transmission conditions. There exists a constant $\alpha_0$ such that for $\alpha$ to be greater than $\alpha_0$
\begin{eqnarray*}
\lim_{k\to\infty}\sum_{l=1}^I\int_{{\Omega}_l}|e_l^k|_\alpha^2 dx =0,
\end{eqnarray*}
\end{theorem}
\begin{proof} For all $l$ in $\{1,I\}$, denote by $\tilde{\Omega}_l$ to be the open set $\Omega_l\backslash\overline{\cup_{l'\in J_l}\Omega_{l'}}$. Let $\varphi_l^k$ be functions in $H^1(\tilde{\Omega}_l)$ and $\varphi_{l}^{k+1}$ be functions in $H^1({\Omega}_l)$ for all $l$ in $I$ such that $\varphi_l^{k+1}=\varphi_{l'}^k$ on $\Gamma_{l,l'}$ for all $l'$ in $J_l$. Now, use $\varphi_l^{k+1}$ and $\varphi_l^k$ as test functions for $(\ref{2e34})$, and take the sum with respect to $l$ in $\{1,I\}$ the integrals $\int_{\tilde{\Omega}_l}\mathfrak{L}_{lR}(\Phi_{l}^{k})\varphi_l^k$ and $\int_{\tilde{\Omega}_l}\mathfrak{L}_{lR}(\Phi_{l}^{k+1})\varphi_l^{k+1}$, then
\begin{eqnarray}\label{2e35}\nonumber
& &-\sum_{l=1}^I\left\{\int_{\tilde{\Omega}_l}\sum_{i,j=1}^na_{i,j}{\partial_i\Phi_l^k}{\partial_j\varphi_l^k}dx+\sum_{i=1}^n\int_{\tilde{\Omega}_l}B_i^l{\partial_i \Phi_l^k}\varphi_l^kdx+\int_{\tilde{\Omega}_l}C^l\Phi_l^k\varphi_l^kdx\right.\\\nonumber
& &-\sum_{l'\in J_l}\int_{\Gamma_{l',l}}p_{l',l}\Phi_l^k\varphi_l^kd\sigma\\
& &\left.+\int_{\tilde{\Omega}_l}\left\{\int_0^{\infty}\left[\left(\frac{\alpha}{2}+c\right)e_l^k-F(u_l^k)+F(u)\right]\exp(-\alpha t)dt\right\}g_l\varphi_l^{k}dx\right\}\\\nonumber
&=&\sum_{l=1}^I\beta_l\left\{\int_{{\Omega}_l}\sum_{i,j=1}^na_{i,j}{\partial_i\Phi_l^{k+1}}{\partial_j\varphi_l^{k+1}}dx+\int_{{\Omega}_l}\sum_{i=1}^nB_i^l{\partial_i \Phi_l^{k+1}}\varphi_l^{k+1}dx\right.\\\nonumber
& &+\int_{{\Omega}_l}C^l\Phi_l^{k+1}\varphi_l^{k+1}dx+\sum_{l'\in J_l}\int_{\Gamma_{l,l'}}p_{l,l'}\Phi_l^{k+1}\varphi_l^{k+1}d\sigma\\\nonumber
& &\left.+\int_{{\Omega}_l}\left\{\int_0^{\infty}\left[\left(\frac{\alpha}{2}+c\right)e_l^{k+1}-F(u_l^{k+1})+F(u)\right]\exp(-\alpha t)dt\right\}g_l\varphi_l^{k+1}dx\right\}.
\end{eqnarray}
\\\h In the above equality, choose $\varphi_l^{k+1}$ to be $\Phi_l^{k+1}$, then there exists $\varphi_l^{k}$, such that $\varphi_l^{k}=\varphi_{l'}^{k+1}$ on $\Gamma_{l,l'}$ for all $l'$ in $J_l$; moreover, 
$$||\varphi_l^{k}||_{H^1({\Omega}_l)}\leq C\sum_{l'\in J_l}||\varphi_{l'}^{k+1}||_{H^1(\Omega_{l'})}\mbox{ and }||\varphi_l^{k}||_{L^2({\Omega}_{l})}\leq C\sum_{l'\in J_l}||\varphi_{l'}^{k+1}||_{L^2(\Omega_{l'})},$$ 
where $C$ is a positive constant.
\\ With these test functions, the right hand side of $(\ref{2e35})$ is greater than or equal to
\begin{eqnarray}\label{2e36}\nonumber
& &\sum_{l=1}^I\beta_l\left\{\int_{{\Omega}_l}\lambda|\nabla\Phi_l^{k+1}|^2dx-\sum_{i=1}^n\int_{{\Omega}_l}||B_i^l||_{L^\infty(\Omega_l)}\left|{\partial_i \Phi_l^{k+1}}\right||\Phi_l^{k+1}|dx\right.\\\nonumber
& &+\frac{\alpha}{4}\int_{{\Omega}_l}|\Phi_l^{k+1}|^2dx+\sum_{l'\in J_l}\int_{\Gamma_{l,l'}}p_{l,l'}|\Phi_l^{k+1}|^2d\sigma\\\nonumber
& &\left.+\int_{{\Omega}_l}\left\{\int_0^{\infty}\left[\left(\frac{\alpha}{2}+c\right)e_l^{k+1}-F(u_l^{k+1})+F(u)\right]\exp(-\alpha t)dt\right\}g_l\varphi_l^{k+1}dx\right\}\\\nonumber
&\geq&\sum_{l=1}^I\beta_l\left\{\int_{{\Omega}_l}\lambda|\nabla\Phi_l^{k+1}|^2dx-\sum_{i=1}^n\int_{{\Omega}_l}||B_i^l||_{L^\infty(\Omega_l)}\left|{\partial_i \Phi_l^{k+1}}\right||\Phi_l^{k+1}|dx\right.\\
& &\left.+\frac{\alpha}{4}\int_{{\Omega}_l}|\Phi_l^{k+1}|^2\right.\\\nonumber
& &\left.+\int_{{\Omega}_l}\left\{\int_0^{\infty}\left[\left(\frac{\alpha}{2}+c\right)e_l^{k+1}-F(u_l^{k+1})+F(u)\right]\exp(-\alpha t)dt\right\}g_l\varphi_l^{k+1}dx\right\}\\\nonumber
&\geq&\sum_{l=1}^I\beta_l\left\{\int_{{\Omega}_l}\frac{\lambda}{2}|\nabla\Phi_l^{k+1}|^2dx+\frac{\alpha}{8}\int_{{\Omega}_l}|\Phi_l^{k+1}|^2\right.\\\nonumber
& &\left.+\int_{{\Omega}_l}\left\{\int_0^{\infty}\left[\left(\frac{\alpha}{2}+c\right)e_l^{k+1}-F(u_l^{k+1})+F(u)\right]\exp(-\alpha t)dt\right\}g_l\varphi_l^{k+1}dx\right\}\\\nonumber
&\geq&\sum_{l=1}^I\beta_l\left\{\int_{{\Omega}_l}\frac{\lambda}{2}|\nabla\Phi_l^{k+1}|^2dx+\frac{\alpha}{8}\int_{{\Omega}_l}|\Phi_l^{k+1}|^2\right.\\\nonumber
& &\left.+\int_{{\Omega}_l}\left[\int_0^{\infty}\left(\frac{\alpha}{2}+c-C'\right)e_l^{k+1}\exp(-\alpha t)dt\right]\left[\int_0^{\infty}e_l^{k+1}\exp(-\alpha t)dt\right]g_l^2dx\right\}\\\nonumber
&\geq&\sum_{l=1}^I\beta_l\left\{\int_{{\Omega}_l}\frac{\lambda}{2}|\nabla\Phi_l^{k+1}|^2dx+\frac{\alpha}{8}\int_{{\Omega}_l}|\Phi_l^{k+1}|^2\right\},
\end{eqnarray}
where 
\begin{eqnarray*}
\left \{ \begin{array}{ll}C'=C \mbox{ if } \int_0^{\infty}e_l^{k+1}\exp(-\alpha t)dt\geq 0,\vspace{.1in}\\
C'=-C\mbox{  if } \int_0^{\infty}e_l^{k+1}\exp(-\alpha t)dt< 0,\end{array}\right. 
\end{eqnarray*}
and notice that $\alpha$ is large enough.
\\ Similarly, we can estimate the left hand side of $(\ref{2e35})$, which is in fact less than or equal to
\begin{eqnarray}\label{2e37}\nonumber
& &\sum_{l=1}^I\left\{\int_{\tilde{\Omega}_l}\Lambda|\nabla\Phi_l^k||\nabla\varphi_l^k|dx+\int_{\tilde{\Omega}_l}2\alpha|\Phi_l^k||\varphi_l^k|dx\right.\\\nonumber
& &\left.+\sum_{i=1}^n\int_{\tilde{\Omega}_l}||B_i^l||_{L^{\infty}(\tilde{\Omega}_l)}\left|{\partial_i \Phi_l^k}\right||\varphi_l^k|dx+\sum_{l'\in J_l}\int_{\Gamma_{l',l}}p_{l',l}|\Phi_l^k||\varphi_l^k|d\sigma\right\}\\\nonumber
&\leq&\sum_{l=1}^IM_1\left[\Lambda \left(||\nabla\Phi_l^k||^2_{L^2(\tilde{\Omega}_l)}+||\nabla\varphi_l^k||^2_{L^2(\tilde{\Omega}_l)}\right)+\frac{\alpha}{2}||\Phi_l^k||^2_{L^2(\tilde{\Omega}_l)}+\frac{\alpha}{2}||\varphi_l^k||^2_{L^2(\tilde{\Omega}_l)}\right.\\
& &+\frac{1}{2}\left(||\nabla\Phi_l^k||^2_{L^2(\tilde{\Omega}_l)}+(\max_{i\in\{1,I\}}||B_i^l||_{L^{\infty}(\tilde{\Omega}_l)})^2||\varphi_l^k||_{L^2(\tilde{\Omega}_l)}^2\right)\\\nonumber
& &\left.+\sum_{l'\in J_l}||p_{l',l}||_{L^\infty(\Gamma_{l',l})}\left(||\Phi_l^k||^2_{H^1(\tilde{\Omega}_l)}+||\varphi_l^k||^2_{H^1(\tilde{\Omega}_l)}\right)\right],
\end{eqnarray}
where $M_1$ is a positive constant depending only on $\{\Omega_l\}_{l\in\{1,I\}}$ and the coefficients of $(\ref{2e3})$. Since $\alpha$ can be chosen such that $\alpha>(\max_{i\in\{1,I\}}||B_i^l||_{L^{\infty}(\tilde{\Omega}_l)})^2$, there exists $M_2$ positive, depending only on $\{\Omega_l\}_{l\in\{1,I\}}$ and the coefficients of $(\ref{2e3})$ such that the right hand side of $(\ref{2e37})$ is less than
\begin{eqnarray}\label{2e38}
& &\sum_{l=1}^IM_2\left[\int_{\tilde{\Omega}_l}\left(\frac{\lambda}{2}|\nabla\Phi_l^{k}|^2dx+\frac{\alpha}{8}|\Phi_l^{k}|^2+\frac{\lambda}{2}|\nabla\Phi_l^{k+1}|^2+\frac{\alpha}{8}|\Phi_l^{k+1}|^2\right)dx\right]\\\nonumber
&\leq&\sum_{l=1}^IM_2\left(\frac{\lambda}{2}||\nabla\Phi_l^{k}||^2_{L^2(\Omega_l)}+\frac{\alpha}{8}||\Phi_l^{k}||^2_{L^2(\Omega_l)}+\frac{\lambda}{2}||\nabla\Phi_l^{k+1}||^2_{L^2(\Omega_l)}+\frac{\alpha}{8}||\Phi_l^{k+1}||^2_{L^2(\Omega_l)}\right).
\end{eqnarray}
\h Define 
\begin{eqnarray}\label{2e39}
E_k:=\sum_{l=1}^I\left(\frac{\lambda}{2}||\nabla\Phi_l^{k}||^2_{L^2(\Omega_l)}+\frac{\alpha}{8}||\Phi_l^{k}||^2_{L^2(\Omega_l)}\right),
\end{eqnarray}
then from $(\ref{2e36})$, $(\ref{2e37})$ and $(\ref{2e38})$, 
\begin{eqnarray}\label{2e40}
(\beta-M_2)E_{k+1}\leq M_2E_k,
\end{eqnarray}
where 
$\beta=\min\{\beta_1,\dots,\beta_I\}$. 
\\ Since $M_2$ depends only on $\{\Omega_l\}_{l\in\{1,I\}}$ and the coefficients of $(\ref{2e3})$, $\beta$ can be chosen large enough, such that $$M_3:=\frac{M_2}{\beta-M_2}<1.$$
We obtain
\begin{eqnarray*}
E_{k}&\leq &M_3^{k}E_0\\\nonumber
& \leq &M_3^{k}\sum_{l=1}^I\left(\frac{\lambda}{2}||\nabla\Phi_l^{0}||^2_{L^2(\Omega_l)}+\frac{\alpha}{8}||\Phi_l^{0}||^2_{L^2(\Omega_l)}\right).
\end{eqnarray*}
That implies
\begin{eqnarray}\label{2e41}
||\Phi_l^{k}||^2_{L^2(\Omega_l)} \leq M_3^{k}\sum_{l=1}^I\left(\frac{4\lambda}{\alpha}||\nabla\Phi_l^{0}||^2_{L^2(\Omega_l)}+||\Phi_l^{0}||^2_{L^2(\Omega_l)}\right).
\end{eqnarray}
Notice that $(\ref{2e41})$ still holds if $M_3$ and $\lambda$ are fixed, and $\alpha$ is replaced by all $y$ which is larger than $\alpha$. This observation leads to
\begin{eqnarray}\label{2e42}
& &\sum_{l=1}^I\int_{{\Omega}_l}\left(\int_0^{\infty}e_l^k\exp(-y t)dtg_l\right)^2dx \\\nonumber
& \leq &M_3^{k}\left[\frac{4\lambda}{y}\sum_{l=1}^I\int_{{\Omega}_l}\left(\int_0^{\infty}|\nabla e_l^0|\exp(-y t)dt\right)^2g_l^2dx\right. \\\nonumber
& &+\frac{4\lambda}{y}\sum_{l=1}^I\int_{{\Omega}_l}\left(\int_0^{\infty}e_l^0\exp(-y t)dt\right)^2|\nabla g_l|^2dx\\\nonumber
& &\left.+\sum_{l=1}^I\int_{{\Omega}_l}\left(\int_0^{\infty}e_l^0\exp(-y t)dt\right)^2g_l^2dx\right].
\end{eqnarray}
Let $\alpha'$ be a constant larger than or equal to $\alpha$, we obtain from the previous inequality that
\begin{eqnarray}\label{2e43}
& &\sum_{l=1}^I\int_{{\Omega}_l}\int_{\alpha'}^{\alpha'+1}\left(\int_0^{\infty}e_l^k\exp(-y t)dt\right)^2g_l^2dydx \\\nonumber
& \leq &M_3^{k}\left[\sum_{l=1}^I\int_{{\Omega}_l}\int_{\alpha'}^{\alpha'+1}\frac{4\lambda}{y}\left(\int_0^{\infty}|\nabla e_l^0|\exp(-y t)dt\right)^2g_l^2dydx\right. \\\nonumber
& &+\sum_{l=1}^I\int_{{\Omega}_l}\int_{\alpha'}^{\alpha'+1}\frac{4\lambda}{y}\left(\int_0^{\infty} e_l^0\exp(-y t)dt\right)^2|\nabla g_l|^2dydx\\\nonumber
& &\left.+\sum_{l=1}^I\int_{{\Omega}_l}\int_{\alpha'}^{\alpha'+1}\left(\int_0^{\infty}e_l^0\exp(-y t)dt\right)^2g_l^2dydx\right].
\end{eqnarray}
\h Using the fact that $u^0$ belongs to $C^{\infty}_c(\overline{\Omega\times(0,\infty)})$, we can infer that the right hand side of $(\ref{2e43})$ is bounded by a constant $M_3^kM_4(\alpha)$. Since $g_l$ is greater than $1$, then
\begin{eqnarray}\label{2e44}
\sum_{l=1}^I\int_{{\Omega}_l}\int_{\alpha'}^{\alpha'+1}\left(\int_0^{\infty}e_l^k\exp(-y t)dt\right)^2dydx \leq M_3^{k}M_4(\alpha).
\end{eqnarray}
\h $(\ref{2e44})$ infers
\begin{eqnarray}\label{2e45}
\lim_{k\to\infty}\sum_{l=1}^I\int_{{\Omega}_l}|e_l^k|_\alpha^2 dx =0,
\end{eqnarray}
that concludes the proof.

\end{proof}

\section{Convergence for Semilinear Elliptic Equations}
\h Consider the semilinear elliptic equation
\begin{equation}
\label{3e1}
\left \{ \begin{array}{ll}-\sum_{i,j=1}^n\partial_j(a_{i,j}{\partial_{i} u})+\sum_{i=1}^nb_{i}{\partial_i u}+cu=F(x,u),\mbox{ in } \Omega,\vspace{.1in}\\ 
u=g,\mbox{ on } \partial\Omega,\end{array}\right. 
\end{equation}
where $\Omega$ is a bounded and smooth domain in  $\mathbb{R}^{n}$. We impose  on the coefficients of $(\ref{3e1})$ Conditions $(A1)$, $(A2)$ in the previous section and the following condition 
\\ $(A3')$ There exists $C>0$, such that  $C<c(x)$ on $\overline{\Omega}$ and for all $x$ in $\mathbb{R}^n$: $$|F(x,z)-F(x,z')|\leq C|z-z'|, \forall z, z'\in \mathbb{R}.$$
 With Conditions $(A1)$, $(A2)$ and $(A3')$, Equation $(\ref{3e1})$ has a unique solution $u$ in $W^{1,2}(\Omega)\cap L^{\infty}(\Omega)$ (see \cite{GilbargTrudinger:2001:EPD}, \cite{Mizohata:1973:TPD}).
\\\h We impose the same way of dividing the domain $\Omega$ and the same notations as in the previous section.
\\ The Schwarz algorithm at the iterate $\#k$ in the $l$-th domain, denoted by $u_{l}^{k}$, is then defined by
\begin{equation}\label{3e2}
\left \{
\begin{array}{ll}
-\sum_{i,j=1}^n\partial_j(a_{i,j}{\partial_{i}  u^k_l})+\sum_{i=1}^nb_{i}{\partial_i  u^k_l}+c u^k_l=F(x,u^k_l),\mbox{ in }\Omega_l,\vspace{.1in}\\
 \mathfrak{B}_{l,l'}u^k_l=\mathfrak{B}_{l,l'}u^{k-1}_{l'}, \mbox{ on }\Gamma_{l,l'},\forall l'\in J_l,\end{array}\right. 
\end{equation}
where $\mathfrak{B}_{l,l'}$ are  either Dirichlet or Robin transmission operators. 
\\ Each iterate also inherits the boundary conditions of $u$:
\[
u_l^k= g \mbox{ on }\partial\Omega_l\cap\partial\Omega.\]
A bounded initial guess $u^0$ in $C^{\infty}(\overline{\Omega\times(0,\infty)})$ is also provided
\begin{eqnarray*}
\mathfrak{B}_{l,l'}u^1_l=u^{0} &\mbox{ on }\Gamma_{l,l'},\forall l'\in J_l.
\end{eqnarray*}
\\\h The difference $e_l^k$ between $u_l^k$ and $u$ is a solution of 
\begin{equation}\label{3e3}
\left \{
\begin{array}{ll}
-\sum_{i,j=1}^n\partial_j(a_{i,j}{\partial_{i} e^k_l})+\sum_{i=1}^nb_{i}{\partial_i  e^k_l}\vspace{.1in}+c e^k_l=F(x,u^k_l)-F(x,u),\mbox{ in }\Omega_l,\vspace{.1in}\\
 \mathfrak{B}_{l,l'}e^k_l=\mathfrak{B}_{l,l'}e^{k-1}_{l'}, \mbox{ on }\Gamma_{l,l'},\forall l'\in J_l.\end{array}\right. 
\end{equation}
Moreover, 
\[
e_l^k= 0 \mbox{ on }\partial\Omega_l.
\]
\subsection{Classical Schwarz Methods}
By induction, each subproblem $(\ref{3e2})$ in each iteration has a unique solution in $W^{1,2}(\Omega)\cap L^{\infty}(\Omega)$ for the Dirichlet transmission condition. The algorithm is well-posed. 
\\ Consider $(\ref{3e3})$ and let $g$ be a bounded and strictly positive function in $C^{2}(\mathbb{R}^n,\mathbb{R})$. Define the following function
\begin{eqnarray*}
\Phi_l^k(x):=(e_l^k(x))^2g(x).
\end{eqnarray*}
Let $c_i$ be $b_i+\sum_{j=1}^n2a_{i,j}{\partial_j g}g^{-1},$ then $c_i$ belongs to $L^{\infty}(\Omega_l),$ and define
\begin{eqnarray}\label{3e9}
\mathfrak{L}_{lD}(\Phi)&=&-\sum_{i,j=1}^n\partial_j(a_{i,j}{\partial_{i}\Phi})+\sum_{i=1}^{n}c_i{\partial_i\Phi}.
\end{eqnarray}
\begin{lemma}\label{linear3} Choose $g$ to be $\tilde{g}^{-1}$ where $\tilde{g}$ is a solution in $C^{2}(\Omega)\cap C(\overline{\Omega})$ of the following equation 
\begin{equation}\label{3e8}
\left \{
\begin{array}{ll}
\sum_{i,j=1}^n\partial_j(a_{i,j}{{\partial_{i} \tilde{g}}})+2(C+||c||_{\infty})\tilde{g}-\sum_{i=1}^{n}b_i{{\partial_i \tilde{g}}}\leq0, \mbox{ in }\Omega_l,\vspace{.1in}\\
 \tilde{g} \mbox{ is strictly positive and bounded on }\overline\Omega,\end{array}\right. 
\end{equation}
then $\mathfrak{L}_{lD}(\Phi_l^k)\leq 0$ in the distributional sense.
\end{lemma}
\begin{proof}
%\\\h  Derive the function $\Phi_l^k$, 
%\begin{eqnarray}\label{3e4}
%{\partial_i\Phi_l^k}& = &2{\partial_i e_l^k}e_l^kg+(e_l^k)^2{\partial_i g}, \forall i\in\{1,\dots,n\},\\\nonumber
%{\partial_{i,j}\Phi_l^k}& = &2{\partial_{i,j} e_l^k}e_l^kg+2{\partial_i e_l^k}{\partial_j e_l^k}g+2{\partial_i e_l^k}e_l^k{\partial_j g}\\\nonumber
%& &+2e_l^k{\partial_j e_l^k}{\partial_i g}+(e_l^k)^2{\partial_{i,j} g},\forall i,j\in\{1,\dots,n\}.
%\end{eqnarray}
Define the operator
\begin{eqnarray*}
\mathfrak{L}_{lD0}:=-\sum_{i,j=1}^n\partial_j(a_{i,j}{\partial_{i}}).
\end{eqnarray*}
A complicated but easy computation gives
\begin{eqnarray}\label{3e5}
\mathfrak{H}(\Phi_l^k)& =&\left(-\sum_{i,j=1}^n2\partial_j(a_{i,j}{\partial_{i}e_l^k})\right)e_l^kg-\sum_{i,j=1}^n2a_{i,j}{\partial_{i} e_l^k}{\partial_j e_l^k}g-\\\nonumber
& &-\sum_{i,j=1}^n4a_{i,j}{\partial_i e_l^k}e_l^k{\partial_j g}+(e_l^k)^2\left(\sum_{i,j=1}^n(-a_{i,j}{\partial_{i,j} g}+\partial_j a_{i,j}\partial_i g)\right).
\end{eqnarray}
\\ That implies
\begin{eqnarray}\label{3e6}
\mathfrak{H}(\Phi_l^k)+\sum_{i=1}^{n}\left(b_i+\sum_{j=1}^n2a_{i,j}{\partial_j g}g^{-1}\right){\partial_i\Phi_l^k}\leq(e_l^k)^2\mathfrak{M},
\end{eqnarray}
where 
\begin{eqnarray}\label{3e7}\nonumber
\mathfrak{M}%&=&-\sum_{i,j=1}^n\partial_j(a_{i,j}{\partial_{i} g})+2g(C+||c||_{\infty})+\sum_{i=1}^{n}b_i{\partial_i g}+\sum_{i,j=1}^n2a_{i,j}{\partial_j g}{\partial_i g}g^{-1}\\
=\left(\sum_{i,j=1}^n\partial_j(a_{i,j}{{\partial_{i} (g^{-1})}})+2(C+||c||_{\infty})(g^{-1})-\sum_{i=1}^{n}b_i{{\partial_i (g^{-1})}}\right)g^2.
\end{eqnarray}
 Therefore, the nonlinear equation $(\ref{3e3})$ has been transformed into the following linearized inequation  of $\Phi_l^k$: 
\begin{eqnarray}\label{3e10}
-\sum_{i,j=1}^n\partial_j(a_{i,j}{\partial_{i}\Phi_l^k})+\sum_{i=1}^{n}c_i(x,t){\partial_i\Phi_l^k}\leq0.
\end{eqnarray}
\end{proof}
\begin{theorem}\label{3t1} Consider the Schwarz algorithm with Dirichlet transmission condition, 
\begin{eqnarray*}
\mathop{\lim}_{k\to\infty}\max_{l\in \{1,\dots,I\}}||u_l^k-u||_{L^{\infty}(\Omega_l)}=0.
\end{eqnarray*}
\end{theorem}
\begin{proof}
{\bf Step 1:} Construct some estimates of the errors $\{e_l^k\}$.
\\\h Consider $(\ref{3e10})$ and define 
\begin{eqnarray}\label{11}
M=\mathop{esssup}_{\partial\Omega_l}\Phi_l^k(x).
\end{eqnarray}
By the weak maximum principle, Theorem $8.1$, \cite{GilbargTrudinger:2001:EPD}, $\Phi_l^k(x)$ is bounded by $M$ almost every where on $\Omega_l$, that implies the following estimates, 
\begin{eqnarray}\label{3e12}
\left(e_l^k\right)^2 g\leq \max_{l'\in J_l}\left(\mathop{esssup}_{\Gamma_{l,l'}}\left(e_{l}^k\right)^2g\right),\mbox{ a.e. in } \Omega_l, \forall l \in \{1,\dots,I\}.
\end{eqnarray}
{\\\bf Step 2:} Convergence of the Algorithm.
\\\h Denote
\begin{eqnarray}\label{3e13}
E^k=\max_{l\in I}\left(\mathop{esssup}_{\Omega_l}\left(e_l^k\right)^2 \right).
\end{eqnarray}
\h From $(\ref{2e19})$, for all $l'$ in $J_l$, and for a.e. $x$ in ${\Gamma_{l,l'}}$
\begin{eqnarray}\label{3e14}
\left(e_{l}^k(x)\right)^2g(x)\leq \max_{l''\in J_{l'}}\left(\mathop{esssup}_{\Gamma_{l',l''}}\left(e_{l'}^{k-1}(x)^2g(x)\right)\right),
\end{eqnarray}
or
\begin{eqnarray}\label{3e15}
\left(e_{l}^k(x)\right)^2& \leq & \frac{1}{g(x)}\max_{l''\in J_{l'}}\left(\mathop{esssup}_{\Gamma_{l',l''}}\left(e_{l'}^{k-1}(x)^2g(x)\right)\right)\\\nonumber
&\leq&\tilde{g}(x)\max_{l''\in J_{l'}}\left(\mathop{esssup}_{\Gamma_{l',l''}}\left(e_{l'}^{k-1}(x)^2\tilde{g}(x)^{-1}\right)\right).
\end{eqnarray}
\\\h Fix $l$ in $\{1,\dots,I\}$ and let $f$ be a function in $C^2(\bar\Omega)$ such that 
%\begin{itemize}
%\item $f>0$ on $\Omega_{l'}$, $\forall l'\in J_l$;
%\item  ${\partial_i f}\ne 0$ on $\overline\Omega_{l'}$, $\forall i\in\{1,\dots,n\}$, $\forall l'\in J_l$;
%\item  $f=0$ on $\partial\Omega_{l'}$, $\forall l'\in J_l$;
%\end{itemize} 
%or
\begin{itemize}
\item $f>0$ on $\Omega_{l'}$, $\forall l'\in J_l$;
\item  There exist two positive real numbers $\epsilon$ small and $M$ large such that for $|\nabla f(x)|<\epsilon$, $\sum_{i,j=1}^n\partial_j(a_{i,j}{{\partial_{i} f}})(x)>M$;
\item  $f=0$ on $\partial\Omega_{l'}$, $\forall l'\in J_l$;
\item $||f||_{\infty}=1$. (We can construct this function by constructing a function $g$ which satisfies the fist three properties, and then take $f={g}\slash {||g||_{\infty}}$).
\end{itemize} 
Let $\rho$ be a constant and put $$\tilde{g}=M_0-\exp(\rho f)=\exp(2\rho)-\exp(\rho f),$$ then, 
\begin{eqnarray}\label{3e16}\nonumber
& &\sum_{i,j=1}^n\partial_j(a_{i,j}{{\partial_{i} \tilde{g}}})+2(C+||c||_{\infty})\tilde{g}-\sum_{i=1}^{n}b_i{{\partial_i \tilde{g}}}\\\nonumber
&=&-\sum_{i,j=1}^n\rho^2\exp(\rho f)\partial_j(a_{i,j}{{\partial_i f}})-\sum_{i,j=1}^na_{i,j}\rho\exp(\rho f){{\partial_{i,j} f}}\\\nonumber
& &+2(C+||c||_{\infty})(\exp(2\rho)-\exp(\rho f))+\sum_{i=1}^{n}b_i\exp(\rho f)\rho {\partial_i f}\\
&<&\exp(-\rho f)\left[-\rho^2M+\epsilon\rho+2(C+||c||_{\infty})\left(\frac{\exp(2\rho)}{\exp(\rho f)}-1\right)+\max_i||b_i||_{\infty}\epsilon\rho\right]\\\nonumber
&<&\exp(-\rho f)\left[-\rho^2M+\epsilon\rho+2(C+||c||_{\infty})\exp(2\rho )+\max_i||b_i||_{\infty}\epsilon\rho\right]\\\nonumber
&<& 0,
\end{eqnarray}
when $\rho$ is large enough and $M>\exp(2\rho)$.
\\ Inequality $(\ref{3e15})$ then becomes
\begin{eqnarray}\label{3e17}
\left(e_{l}^k(x)\right)^2& \leq & (M_0-\exp(\rho f(x)))\max_{l''\in J_{l'}}\left(\mathop{esssup}_{\Gamma_{l',l''}}\left(\left(e_{l'}^{k-1}(x)\right)^2(M_0-1)^{-1}\right)\right).
\end{eqnarray}
Since $\Gamma_{l,l'}$ lies inside $\Omega_{l'}$, $f(x)$ is strictly positive on $\Gamma_{l,l'}$. Hence, there exists $M_{1,l}$ strictly less than $1$, such that  
\begin{eqnarray}\label{3e18}
\left(e_{l}^k(x)\right)^2\leq M_{1,l}\max_{l''\in J_{l'}}\left(\mathop{esssup}_{\Gamma_{l',l''}}\left(e_{l'}^{k-1}(x)\right)^2\right)\leq M_{1,l}E^{k-1},\forall l' \in J_l,  \mbox{ for a.e  $x$ in }{\Gamma_{l,l'}}.
\end{eqnarray}
\\\h Similarly as in $(\ref{3e16})$, let $f$ be a function in $C^2(\bar\Omega)$ such that 
%\begin{itemize}
%\item $f>0$ on $\Omega_{l}$; 
%\item $\forall i\in\{1,\dots,n\}$, ${\partial_i f}\ne 0$ on $\overline\Omega_{l}$; 
%\item $f=0$, on $\partial\Omega_{l}$; 
%\end{itemize}
%or
\begin{itemize}
\item  $f=0$ on $\partial\Omega_{l}$;
\item $f>0$ on on $\Omega_{l}$;
\item  There exist two positive constants $\epsilon$ small enough and $M$ large enough such that for $|\nabla f_l(x)|<\epsilon$, we have that $\sum_{i,j=1}^n\partial_j(a_{i,j}{{\partial_{i} f}})(x)>M$;
\item $||f||_{\infty}=1$.
\end{itemize} 
and let $\rho$ be a constant large enough. Setting $$\tilde{g}=M_0-\exp(\rho f),$$ 
\begin{eqnarray}\label{3e19}\nonumber
\left(e_l^k(x)\right)^2 &\leq& \tilde{g}(x)\max_{l'\in J_l}\left(\sup_{\Gamma_{l,l'}}\left(e_{l}^k(x)\right)^2(M_0-1)\right)\\
&\leq& \max_{l'\in J_l}\left(\sup_{\Gamma_{l,l'}}\left(e_{l}^k(x)\right)^2\right),\mbox{ for a.e  $x$ in }\Omega_l.
\end{eqnarray}
\\\h Combining $(\ref{3e18})$ and $(\ref{3e19})$, we can deduce that there exists $M_2$ strictly less than $1$  satisfying
\begin{eqnarray}\label{3e20}
E^k\leq M_2E^{k-1},
\end{eqnarray}
that leads to $$\lim_{k\to\infty}E^k=0.$$ 
\end{proof}
\subsection{Optimized Schwarz Methods}
The optimized Schwarz waveform relaxation algorithms are defined by replacing the Dirichlet transmission operators by Robin ones
$$\mathfrak{B}_{l,l'}v=\sum_{i,j=1}^n a_{i,j}{\partial_i v}n_{l,l',j}+p_{l,l'}v,$$ 
where $n_{l,l',j}$ is the $j$-th component of the outward unit normal vector of $\Gamma_{l,l'}$; $p_{l,l'}$ is positive and belongs to $L^{\infty}(\Gamma_{l,l'})$.
By induction argument, each subproblem $(\ref{2e2})$ in each iteration has a unique solution in $H^1(\Omega)$ and the algorithm is well-posed. 
\\\h { In general, optimized Schwarz algorithms do not always converge when  applied to elliptic equations, as shown in the following example.} 
\begin{example}\label{ex1}
Consider the elliptic problem on the domain $\Omega=(0,L)$
\begin{equation}\label{3eex1}
\left \{
\begin{array}{ll}
u''-3u'-4u=f,&\mbox{ in }\Omega,\vspace{.1in}\\
 u(0)=u(L)=0,\end{array}\right. 
\end{equation}
where $f$ belongs to $C^{\infty}([0,L]).$ Divide $\Omega$ into two subdomains $\Omega_1=(0,L_2)$ and $\Omega_2=(L_1,L)$, with $0<L_1<L_2<L$, and consider the domain decomposition algorithm
\begin{eqnarray*}
\left \{
\begin{array}{ll}
(u_1^{k+1})''-3(u_1^{k+1})'-4u_1^{k+1}=f,\mbox{ in }(0,L_2),\vspace{.1in}\\
 u_1^{k+1}(0)=0 \mbox{ and } (u_1^{k+1})'(L_2)+pu_1^{k+1}(L_2)=(u_2^{k})'(L_2)+pu_2^{k}(L_2),\end{array}\right. 
\end{eqnarray*}
\begin{eqnarray*}
\left \{
\begin{array}{ll}
(u_2^{k+1})''-3(u_2^{k+1})'-4u_2^{k+1}=f,\mbox{ in }(L_1,L),\vspace{.1in}\\
 u_2^{k+1}(L)=0 \mbox{ and } (u_2^{k+1})'(L_1)-qu_2^{k+1}(L_1)=(u_1^{k})'(L_1)-qu_1^{k}(L_1),\end{array}\right. 
\end{eqnarray*}
where $p$, $q$ are positive numbers.
\\ The errors $e_1^k$ and $e_2^k$ from the above equations can be obtained $$e_1^{k+1}=A_{k+1}(\exp(4x)-\exp(-x)),$$
$$e_2^{k+1}=B_{k+1}(\exp(4(x-L))-\exp(-(x-L))),$$
where
\begin{eqnarray*}
\tau_1=\frac{A_{k+1}}{B_k}=\frac{4\exp(4(L_2-L))+\exp(-(L_2-L))+p(\exp(4(L_2-L))-\exp(-(L_2-L)))}{4\exp(4L_2)+\exp(-L_2)+p(\exp(4L_2)-\exp(-L_2))},\\
\tau_2=\frac{B_{k+1}}{A_k}=\frac{4\exp(4L_1)+\exp(-L_1)-q(\exp(4L_1)-\exp(-L_1))}{4\exp(4(L_1-L))+\exp(-(L_1-L))-q(\exp(4(L_1-L))-\exp(-(L_1-L)))}.
\end{eqnarray*}
Set $$\tau=\left|\frac{A_{k+1}B_{k+1}}{B_kA_k}\right|,$$
then
\begin{eqnarray}\label{3eex2}
\tau&=&\left|\frac{4\exp(5L_2)+\exp(5L)+p(\exp(5L_2)-\exp(5L))}{4\exp(5L_2)+1+p(\exp(5L_2)-1)}\right|\\\nonumber
& &\times\left|\frac{4\exp(5L_1)+1-q(\exp(5L_1)-1)}{4\exp(5L_1)+\exp(5L)-q(\exp(5L_1)-\exp(5L))}\right|.
\end{eqnarray}
The algorithm converges if and only if  $\tau$  is smaller than $1$.
\\\h For $p=1$ and $q$ large, 
\begin{eqnarray}\label{3eex3}
\tau\backsim\frac{\exp(5L_1)-1}{-\exp(5L_1)+\exp(5L)}.
\end{eqnarray}
Since the right hand side of $(\ref{3eex3})$ is greater than $1$ for $L_1>L\slash\ln2$ and $L$ large, the algorithm does not converge for $p=1$ and $q$ large.
\\ Note that 
$$\tau_1=\left|\frac{4\exp(5L_2)+\exp(5L)+p(\exp(5L_2)-\exp(5L))}{4\exp(5L_2)+1+p(\exp(5L_2)-1)}\right|,$$
can be made larger than 1 since 
$$\left|\frac{4\exp(5L_2)+\exp(5L)}{4\exp(5L_2)+1}\right|$$
is larger than 1;
\\ and the second term can be made larger than 1 since
$$\tau_2=\frac{\exp(5L_1)-1}{-\exp(5L_1)+\exp(5L)}$$
is larger than $1$ for $L_1>L\slash\ln2$ and $L$ large.
\end{example}
\begin{remark}
In the above example, the Schwarz algorithms converge if $p$ and $q$ are large. This observation leads to Theorem $\ref{3t2}$ below. Moreover, the convergence factor is still small if $p$, $q$ are chosen to be negative. At least, the behavior of $\tau$ is quite the same when $q$ tend to $+\infty$ and $-\infty$. This observation is different from what was seen from the convergence rate in the case where the subdomains are two half lines in \cite{Bennequin:2009:AHB} and \cite{Gander:2007:OSW}.
\end{remark}
Introduce the modified Robin transmission operator, based on $(T_2)$
$$\mathfrak{B}^\rho_{l,l'}v=\sum_{i,j=1}^n a_{i,j}{\partial_i v}n_{l,l',j}+\rho p_{l,l'}v,$$ 
where $\rho$ is a positive parameter.
\begin{theorem}\label{3t2} Consider Schwarz Algorithms with Robin transmission conditions, if we replace $\mathfrak{B}_{l,l'}$ by $\mathfrak{B}^\rho_{l,l'}$, then there exists $\rho_0$ such that for $\rho>\rho_0$, the algorithms converge in the following sense
\begin{eqnarray*}
\mathop{\lim}_{k\to\infty}\max_{l\in\{1,\dots, I\}}||u_l^k-u||_{L^{2}(\Omega_l)}=0.
\end{eqnarray*}
\end{theorem}
\begin{remark} Consider again Example 3.1, and a Schwarz algorithm which diverges with $p=0$ and $q=q_0$, there exists $L_1$, $L_2$, $L$ such that the algorithm diverges
\begin{eqnarray*}
\left \{
\begin{array}{ll}
(u_1^k)''-3(u_1^k)'-4u_1^k=f,\mbox{ in }(0,L_2),\vspace{.1in}\\
 u_1^k(0)=0 \mbox{ and } (u_1^k)'(L_2)=(u_2^{k-1})'(L_2),\end{array}\right. 
\end{eqnarray*}
\begin{eqnarray*}
\left \{
\begin{array}{ll}
(u_2^k)''-3(u_2^k)'-4u_2^k=f,\mbox{ in }(L_1,L),\vspace{.1in}\\
 u_2^k(L)=0 \mbox{ and } (u_2^{k})'(L_1)-q_0u_2^{k}(L_1)=(u_1^{k-1})'(L_1)-q_0u_1^{k-1}(L_1),\end{array}\right. 
\end{eqnarray*}
 where $q_0$ is a large constant.\\
Let $P$ be a function in $C^1(\mathbb{R})$, and put $w=u\exp(P)$, Equation $(\ref{3eex1})$ can be transformed into
\begin{equation}\label{3eex5}
\left \{
\begin{array}{ll}
w''-(3+2P')w'+(-4-3P'+(P')^2-P'')w=f,&\mbox{ in }(0,L),\vspace{.1in}\\
 w(0)=w(L)=0.\end{array}\right. 
\end{equation} 
The Schwarz algorithm then becomes
\begin{eqnarray*}
\left \{
\begin{array}{ll}
(w_1^k)''-(3+2P')(w_1^k)'+(-4-3P'+(P')^2-P'')w_1^k=f,\mbox{ in }(0,L_2),\vspace{.1in}\\
 w_1^k(0)=0 \mbox{ and } ((w_1^k)'-P'w_1^k)(L_2)=((w_2^{k-1})'-P'w_2^{k-1})(L_2),\end{array}\right. 
\end{eqnarray*}
\begin{eqnarray*}
\left \{
\begin{array}{ll}
(w_2^k)''-(3+2P')(w_2^k)'+(-4-3P'+(P')^2-P'')w_2^k=f,\mbox{ in }(L_1,L),\vspace{.1in}\\
 w_2^k(L)=0 \mbox{ and } (w_2^k)'(L_1)-(P'+q_0)w_2^k(L_1)=((w_1^{k-1})'(L_1)-(P'+q_0)w_2^k(L_1),\end{array}\right. 
\end{eqnarray*}
 We can deduce that given a pair of numbers $(p,q)$, we can find a class of  functions $P$ such that $-P'(L_2)=p$ and $-P'(L_1)=q+q_0$, and the Schwarz algorithm with the associated equation $(\ref{3eex5})$ and this Robin transmission condition does not converge. However, Theorem $\ref{3t1}$ announces that we can make the algorithms converge, even if they do not converge initially, by increasing the parameter $\rho$.
\end{remark}
\begin{proof}
{\bf Step 1:} Linearize the equation $(\ref{2e3})$.
\\\h Consider the equation $(\ref{2e3})$ and let $g_l$ be a strictly positive bounded function in $C^2(\Omega_l,\mathbb{R})$. Define the following function
\begin{eqnarray*}
\Phi_l^k(x):=e_l^k(x)g_l(x).
\end{eqnarray*}
A complicated but easy computation gives
\begin{eqnarray}\label{3e26}
0&=&-\sum_{i,j=1}^n\partial_j(a_{i,j}{\partial_{i} \Phi_l^k})+\sum_{i=1}^nb_{i}{\partial_i \Phi_l^k}+\sum_{i,j=1}^na_{i,j}\left({\partial_i \Phi_l^k}\frac{{\partial_j g_l^k}}{g_l^k}+{\partial_j \Phi_l^k}\frac{{\partial_j g_l^k}}{g_l^k}\right)\\\nonumber
& &+\left(\sum_{i,j=1}^na_{i,j}\frac{{\partial_{i,j} g_l^k}}{g_l^k}-\sum_{i,j=1}^na_{i,j}\frac{2{\partial_i g_l^k}{\partial_j g_l^k}}{(g_l^k)^2}-\sum_{i=1}^nb_{i}\frac{{\partial_i g_l^k}}{g_l^k}+(c-\overline{F})(g_l^k)^{-1}\right)\Phi_l^k,
\end{eqnarray}
where 
\begin{eqnarray*}
\left \{ \begin{array}{ll}\bar{F}(x)=0 \mbox{ if } u_l^k(x)=u(x), x\in\Omega,\vspace{.1in}\\
\bar{F}(x)=\frac{F(u_l^k(x))-F(u(x))}{u_l^k(x)-u(x)}\mbox{  if } u_l^k(x)\ne u(x), x\in\Omega,\end{array}\right. 
\end{eqnarray*}
$\bar{F}$ is then bounded as $F$ is Lipschitz.
\\ Similar as in $(\ref{3e7})$ and in Step 2 of the proof of Theorem $\ref{3t1}$, we rewrite the last term on the right hand side of $(\ref{3e26})$ into the following form
\begin{eqnarray}\label{3e27}
\left(-\sum_{i,j=1}^n\partial_i(a_{i,j}{{\partial_{i} ((g_l^k)^{-1})}})+(c-\overline{F})((g_l^k)^{-1})+\sum_{i=1}^{n}b_i{{\partial_i ((g_l^k)^{-1})}}\right)\Phi_l^k,
\end{eqnarray}
\\ and use the same argument as in $(\ref{3e8})$: choose $g_l^k$ to be $\tilde{g_l}^{-1}$ where $\tilde{g_l}$ is a solution in $C^{2}(\Omega)\cap C(\overline{\Omega})$ of the following equation 
\begin{equation}\label{3e28}
\left \{
\begin{array}{ll}
-\sum_{i,j=1}^n\partial_i(a_{i,j}{{\partial_{i} \tilde{g_l}}})-K\tilde{g_l}+\sum_{i=1}^{n}b_i{{\partial_i \tilde{g_l}}}\geq0, \mbox{ in }\Omega_l,\vspace{.1in}\\
 \tilde{g} \mbox{ is strictly positive and bounded on }\overline\Omega,\end{array}\right. 
\end{equation}
where $K$ is a large enough constant.
\\\h Since $p_{l,l'}$ is strictly positive for all $l$ in $\{1,\dots,I\}$ and $l'$ in $J_l$, there exist functions $f_l$, $l\in\{1,\dots,I\}$, in $C^2(\bar\Omega)$ such that
\begin{itemize}
\item $\sum_{i,j=1}^na_{i,j}n_{l,l',j}{\partial_i f_l}=p_{l,l'}$ on $\Gamma_{l,l'}$, $\forall l'\in J_l$.
\item $\sum_{i,j=1}^na_{i,j}n_{l',l,j}{\partial_i f_l}=p_{l',l}$ on $\Gamma_{l',l}$, $\forall l'$: $l\in J_{l'}$. 
\item  There exist two positive constants $\epsilon$ small enough and $M$ large enough such that for $|\nabla f_l(x)|<\epsilon$, $\sum_{i,j=1}^n\partial_j(a_{i,j}{{\partial_{i} f}})(x)>M$.
\item $f_l=0$, on $\Gamma_{l,l'}$, $\forall l'\in J_l$; and $f_{l'}=\alpha_l$, on $\Gamma_{l',l}$, $\forall l':$ $l\in J_{l'}$. 
\item $||f||_{\infty}=1$. (We can construct this function by constructing a function $g$ which satisfies the fist properties, and then take $f={g}\slash {||g||_{\infty}}$).
\end{itemize}
Similar as in $(\ref{3e16})$, let $\rho$ be a constant large enough and put $\tilde{g}=M_3-\exp(-\rho f)$, where $M_3$ is a positive constant, 
\begin{eqnarray}\label{3e16}\nonumber
& &\sum_{i,j=1}^n\partial_j(a_{i,j}{{\partial_{i} \tilde{g}}})+2(C+||c||_{\infty})\tilde{g}-\sum_{i=1}^{n}b_i{{\partial_i \tilde{g}}}\\\nonumber
&=&-\sum_{i,j=1}^n\rho^2\exp(\rho f)\partial_j(a_{i,j}{{\partial_i f}})-\sum_{i,j=1}^na_{i,j}\rho\exp(\rho f){{\partial_{i,j} f}}\\\nonumber
& &+2(C+||c||_{\infty})(M_3-\exp(\rho f))+\sum_{i=1}^{n}b_i\exp(\rho f)\rho {\partial_i f}\\\nonumber
&=&\left(-\lambda\rho^2M-\sum_{i,j=1}^na_{i,j}{{\partial_{i,j} f}}\rho\right.\\\nonumber
& &\left.+2(C+||c||_{\infty})\frac{M_3-\exp(\rho f)}{\exp(\rho f)}+\sum_{i=1}^{n}b_i\rho {\partial_i f}\right)\exp(\rho f)\\\nonumber
&<&0,
\end{eqnarray}
when $\rho$ is large enough. 
\\\h Denote the right hand side of the equation $(\ref{3e26})$ by $\mathfrak{L}_l(\Phi_l^k)$, then it can be rewritten in the following form
\begin{eqnarray}\label{3e29}
\mathfrak{L}_l(\Phi_l^k)&=&-\sum_{i,j=1}^n\partial_j(a_{i,j}{\partial_{i} \Phi_l^k})+\sum_{i=1}^nB^l_{i}{\partial_i \Phi_l^k}+C^l\Phi_l^k,
\end{eqnarray}
where $B^l_{i}$ and $C^l$ are functions in $L^{\infty}(\mathbb{R}^n)$, $C^l$ is bounded from below by $ K+c+\overline{F}$. $\rho$ can be chosen such that there exists $\alpha$ large enough, $2\alpha>C_l>{\alpha}$.
\\\h Now, consider the Robin transmission condition on the boundary $\Gamma_{l,l'}$
\begin{eqnarray}\label{3e30}
\mathfrak{B}_l(\Phi_l^k)&=&\sum_{i,j=1}^na_{i,j}{\partial_i \Phi_l^k}n_{l,l',j}\\\nonumber
&=&\left(\sum_{i,j=1}^na_{i,j}n_{l,l',j}{\partial_i e_l^k}\right)g_l+\sum_{i,j=1}^na_{i,j}n_{l,l',j}{\partial_i g_l}e_l^k\\\nonumber
&=&\left(\sum_{i,j=1}^na_{i,j}n_{l',l,j}{\partial_i e_l^k}+p_{l,l'}e_l^k\right)g_l\\\nonumber
&=&\left(\sum_{i,j=1}^na_{i,j}n_{l',l,j}{\partial_i e_{l'}^{k-1}}+p_{l,l'}e_{l'}^k\right)g_l\\\nonumber
&=&\left(\sum_{i,j=1}^na_{i,j}{\partial_i \Phi_{l'}^{k-1}}n_{l',l,j}\right)\frac{g_{l}}{g_l'}=\mathfrak{B}_l(\Phi_{l'}^{k-1})\frac{g_{l}}{g_l'}.
\end{eqnarray}
\\ We can choose $f_l$ such that $$\frac{g_{l'}}{g_l}=\beta_l, \mbox{ on } \Gamma_{l,l'}, \forall l'\in J_l,$$ where $\beta_l$ is a constant greater than $1$.
\\\h From the previous calculation on $\mathfrak{B}_l(\Phi_l^k)$ and $\mathfrak{L}_l(\Phi_l^k)$, $\Phi_l^k$ is in fact a solution of the following equation
\begin{equation}\label{3e31}\
\left \{
\begin{array}{ll}
\mathfrak{L}_l(\Phi_l^k)=0,&\mbox{ in }\Omega_l\times(0,\infty),\vspace{.1in}\\
 \beta_l\mathfrak{B}_{l,l'}(\Phi_l^k)=\mathfrak{B}_{l,l'}(\Phi_{l'}^{k-1}) &\mbox{ on }\Gamma_{l,l'}\times(0,\infty),\forall l'\in J_l.\end{array}\right. 
\end{equation}
{\bf Step 2:} The Proof of Convergence.
\\\h Denote by $\tilde{\Omega}_l$ to be the open set $\Omega_l\backslash\overline{\cup_{l'\in J_l}\Omega_{l'}}$. For each $l$ in $\{1,\dots,I\}$, let $\varphi_l^k$ to be a function in $H^1({\Omega}_l)$ and $\varphi_{l}^{k+1}$ to be a function in $H^1({\Omega}_l)$  such that $\varphi_l^{k+1}=\varphi_{l'}^k$ on $\Gamma_{l,l'}$ for all $l'$ in $J_l$. Now, using $\varphi_l^{k+1}$ and $\varphi_l^k$ as test functions for all subdomains, we obtain
\begin{eqnarray}\label{3e32}\nonumber
& &-\sum_{l=1}^I\left\{\int_{\tilde{\Omega}_l}\sum_{i,j=1}^na_{i,j}{\partial_i\Phi_l^k}{\partial_j\varphi_l^k}dx+\int_{\tilde{\Omega}_l}\sum_{i=1}^nB_i^l{\partial_i \Phi_l^k}\varphi_l^kdx+\int_{\tilde{\Omega}_l}C^l\Phi_l^k\varphi_l^kdx\right.\\\nonumber
& &\left.-\sum_{l'\in J_l}\int_{\Gamma_{l',l}}p_{l',l}\Phi_l^k\varphi_l^kd\sigma\right\}\\\nonumber
&=&\sum_{l=1}^I\beta_l\left\{\int_{{\Omega}_l}\sum_{i,j=1}^na_{i,j}{\partial_i\Phi_l^{k+1}}{\partial_j\varphi_l^{k+1}}dx+\int_{{\Omega}_l}\sum_{i=1}^nB_i^l{\partial_i \Phi_l^{k+1}}\varphi_l^{k+1}dx\right.\\
& &\left.+\int_{{\Omega}_l}C^l\Phi_l^{k+1}\varphi_l^{k+1}dx+\sum_{l'\in J_l}\int_{\Gamma_{l,l'}}p_{l,l'}\Phi_l^{k+1}\varphi_l^{k+1}d\sigma\right\}.
\end{eqnarray}
\\\h In the above equality, choose $\varphi_l^{k+1}$ to be $\Phi_l^{k+1}$, then there exists $\varphi_l^{k}$ such that $\varphi_l^{k}=\varphi_{l'}^{k+1}$ on $\Gamma_{l,l'}$ for all $l'$ in $J_l$; and $$||\varphi_l^{k}||_{H^1({\Omega}_l)}\leq C\sum_{l'\in J_l}||\varphi_{l'}^{k+1}||_{H^1(\Omega_{l'})};||\varphi_l^{k}||_{L^2({\Omega}_{l})}\leq C\sum_{l'\in J_l}||\varphi_{l'}^{k+1}||_{L^2(\Omega_{l'})}.$$ With these test functions, the right hand side of $(\ref{3e32})$ is greater than or equal to
\begin{eqnarray}\label{3e33}\nonumber
& &\sum_{l=1}^I\beta_l\left\{\int_{{\Omega}_l}\lambda|\nabla\Phi_l^{k+1}|^2dx-\sum_{i=1}^n\int_{{\Omega}_l}||B_i^l||_{L^\infty(\Omega_l)}\left|{\partial_i \Phi_l^{k+1}}\right||\Phi_l^{k+1}|dx\right.\\\nonumber
& &\left.+\alpha\int_{{\Omega}_l}|\Phi_l^{k+1}|^2dx+\sum_{l'\in J_l}\int_{\Gamma_{l',l}}p_{l',l}|\Phi_l^{k+1}|^2d\sigma\right\}\\\nonumber
&\geq&\sum_{l=1}^I\beta_l\left[\int_{{\Omega}_l}\lambda|\nabla\Phi_l^{k+1}|^2dx-\sum_{i=1}^n\int_{{\Omega}_l}||B_i^l||_{L^\infty(\Omega_l)}\left|{\partial_i \Phi_l^{k+1}}\right||\Phi_l^{k+1}|dx\right.\\\nonumber
& &\left.+\alpha\int_{{\Omega}_l}|\Phi_l^{k+1}|^2\right]\\\nonumber
&\geq&\sum_{l=1}^I\beta_l\left[\int_{{\Omega}_l}\frac{\lambda}{2}|\nabla\Phi_l^{k+1}|^2dx+\frac{\alpha}{2}\int_{{\Omega}_l}|\Phi_l^{k+1}|^2\right],
\end{eqnarray}
with $\alpha$ being large enough.
\\ Similarly, we estimate the left hand side of $(\ref{3e32})$, which is in fact bounded by
\begin{eqnarray}\label{3e34}\nonumber
& &\sum_{l=1}^I\left\{\int_{\tilde{\Omega}_l}\Lambda|\nabla\Phi_l^k||\nabla\varphi_l^k|dx+\int_{\tilde{\Omega}_l}2\alpha|\Phi_l^k||\varphi_l^k|dx\right.\\\nonumber
& &\left.+\sum_{i=1}^n\int_{\tilde{\Omega}_l}||B_i^l||_{L^{\infty}(\tilde{\Omega}_l)}\left|{\partial_i \Phi_l^k}\right||\varphi_l^k|dx+\sum_{l'\in J_l}\int_{\Gamma_{l',l}}p_{l',l}|\Phi_l^k||\varphi_l^k|d\sigma\right\}\\\nonumber
&\leq&\sum_{l=1}^IM_4\left[\Lambda \left(||\nabla\Phi_l^k||^2_{L^2(\tilde{\Omega}_l)}+||\nabla\varphi_l^k||^2_{L^2(\tilde{\Omega}_l)}\right)+{\alpha}||\Phi_l^k||^2_{L^2(\tilde{\Omega}_l)}+{\alpha}||\varphi_l^k||^2_{L^2(\tilde{\Omega}_l)}\right.\\\nonumber
& &+\left(||\nabla\Phi_l^k||^2_{L^2(\tilde{\Omega}_l)}+(\max_{i\in\{1,I\}}||B_i^l||_{L^{\infty}(\tilde{\Omega}_l)})^2||\varphi_l^k||_{L^2(\tilde{\Omega}_l)}^2\right)\\
& &\left.+\sum_{l'\in J_l}||p_{l',l}||_{L^\infty(\Gamma_{l',l})}\left(||\Phi_l^k||^2_{H^1(\tilde{\Omega}_l)}+||\varphi_l^k||^2_{H^1(\tilde{\Omega}_l)}\right)\right],
\end{eqnarray}
where $M_4$ is a positive constant which depends only on $\{\Omega_l\}_{l\in\{1,I\}}$ and the coefficients of $(\ref{3e3})$. Since $\alpha$ can be chosen such that $\alpha>(\max_{i\in\{1,\dots,I\}}||B_i^l||_{L^{\infty}(\tilde{\Omega}_l)})^2$, there exists $M_5$ positive, depending only on $\{\Omega_l\}_{l\in\{1,\dots,I\}}$ and the coefficients of $(\ref{3e3})$ such that the right hand side of $(\ref{3e34})$ is less than
\begin{eqnarray}\label{3e35}
& &\sum_{l=1}^IM_5\left[\int_{\tilde{\Omega}_l}\left(\frac{\lambda}{2}|\nabla\Phi_l^{k}|^2+\frac{\alpha}{2}|\Phi_l^{k}|^2+\frac{\lambda}{2}|\nabla\Phi_l^{k+1}|^2+\frac{\alpha}{2}|\Phi_l^{k+1}|^2\right)dx\right]\\\nonumber
&\leq&\sum_{l=1}^IM_5\left(\frac{\lambda}{2}||\nabla\Phi_l^{k}||^2_{L^2(\Omega_l)}+\frac{\alpha}{2}||\Phi_l^{k}||^2_{L^2(\Omega_l)}+\frac{\lambda}{2}||\nabla\Phi_l^{k+1}||^2_{L^2(\Omega_l)}+\frac{\alpha}{2}||\Phi_l^{k+1}||^2_{L^2(\Omega_l)}\right).
\end{eqnarray}
Define 
\begin{eqnarray}\label{3e36}
E_k:=\sum_{l=1}^I\left(\frac{\lambda}{2}||\nabla\Phi_l^{k}||^2_{L^2(\Omega_l)}+\frac{\alpha}{2}||\Phi_l^{k}||^2_{L^2(\Omega_l)}\right),
\end{eqnarray}
then from $(\ref{3e34})$, $(\ref{3e35})$ and $(\ref{3e36})$, 
\begin{eqnarray}\label{3e37}
(\beta-M_5)E_{k+1}\leq M_5E_k,
\end{eqnarray}
$\beta=\min\{\beta_1,\dots,\beta_I\}$. Since $M_5$ depends only on $\{\Omega_l\}_{l\in{1,I}}$ and the coefficients of $(\ref{2e3})$, $\beta$ can be chosen large enough, such that $\frac{M_5}{\beta-M_5}<1$, then
\begin{eqnarray}\label{3e38}
E_{k}\leq \left(\frac{M_5}{\beta-M_5}\right)^{k-1}E_1,
\end{eqnarray}
which means $E_k$ tends to $0$ as $k$ tends to infinity. This concludes the proof.
\end{proof}
 
\section{Conclusions}
We have introduced a new class of techniques to study the convergene of Schwarz methods. In particular, classical Schwarz methods are proved to converge when being applied to both parabolic and elliptic equations. On the contrary, Schwarz methods with Robin transmission conditions only converge when we use them for parabolic equations, though they were proved to converge faster than classical ones in previous studies. For elliptic equations, we have given a counter example, where we can see that optimized Schwarz methods do not converge; and for each optimized Schwarz algorithm, there exists a class of elliptic equations which is not stable with this algorithm. A new way of stabilizing the algorithms has then been proposed. 	
\\ {\bf Acknowledgements.} The author would like to express his gratitude to his thesis advisor, Professor Laurence Halpern for her very kind help and support. He is grateful to Professor Martin Gander and Dr J\'er\'emie Szeftel for valuable and fruitful discussions on the subject. He would also like to thank Professor Hatem Zaag for valuable advice. 
\bibliographystyle{plain}\bibliography{ProofofConvergenceWeakFormSolution}
\end{document}